# TWO-PLAYER NONZERO–SUM STOPPING GAMES IN DISCRETE TIME

By Eran Shmaya and Eilon Solan[1]

*Tel Aviv University*

We prove that every two-player nonzero–sum stopping game in discrete time admits an $\varepsilon$-equilibrium in randomized strategies for every $\varepsilon > 0$. We use a stochastic variation of Ramsey's theorem, which enables us to reduce the problem to that of studying properties of $\varepsilon$-equilibria in a simple class of stochastic games with finite state space.

**1. Introduction.** The following optimization problem was presented by Dynkin (1969). Two players observe a realization of two real-valued processes $(x_n)$ and $(R_n)$. Player 1 can *stop* whenever $x_n \geq 0$, and player 2 can *stop* whenever $x_n < 0$. At the first stage $\theta$ in which one of the players stops, player 2 pays player 1 the amount $R_\theta$ and the process terminates. If no player ever stops, player 2 does not pay anything.

A *strategy* of player 1 is a stopping time $\mu$ that satisfies $\{\mu = n\} \subseteq \{x_n \geq 0\}$ for every $n \geq 0$. A strategy $\nu$ of player 2 is defined analogously. The termination stage is simply $\theta = \min\{\mu, \nu\}$. For a given pair $(\mu, \nu)$ of strategies, denote by

$$\gamma(\mu, \nu) = \mathbf{E}[\mathbf{1}_{\{\theta < \infty\}} R_\theta]$$

the expected payoff to player 1.

Dynkin (1969) proved that if $\sup_{n \geq 0} |R_n| \in L_1$, this problem has a value $v$; that is,

$$v = \sup_\mu \inf_\nu \gamma(\mu, \nu) = \inf_\nu \sup_\mu \gamma(\mu, \nu).$$

He moreover characterized $\varepsilon$-optimal strategies; that is, strategies $\mu$ (resp. $\nu$) that achieve the supremum (resp. the infimum) up to $\varepsilon$.

Received May 2002; revised February 2004.
[1]Supported by the Israel Science Foundation Grant 69/01-1.
*AMS 2000 subject classifications.* Primary 60G40; secondary 91A15, 91A05.
*Key words and phrases.* Stopping games, Dynkin games, stochastic games, $\varepsilon$-equilibrium, randomized stopping times.

Neveu (1975) generalized this problem by allowing both players to stop at every stage, and by introducing three real-valued processes $(R_{\{1\},n})$, $(R_{\{2\},n})$ and $(R_{\{1,2\},n})$. The expected payoff to player 1 is defined by

$$\gamma(\mu,\nu) = \mathbf{E}[\mathbf{1}_{\{\mu<\nu\}} R_{\{1\},\mu} + \mathbf{1}_{\{\mu>\nu\}} R_{\{2\},\nu} + \mathbf{1}_{\{\mu=\nu<\infty\}} R_{\{1,2\},\mu}].$$

Neveu (1975) then proved that this problem has a value, provided

(a) $\sup_{n\geq 0} \max\{|R_{\{1\},n}|, |R_{\{2\},n}|, |R_{\{1,2\},n}|\} \in L_1$ and
(b) $R_{\{1\},n} = R_{\{1,2\},n} \leq R_{\{2\},n}$.

Recently Rosenberg, Solan and Vieille (2001) studied games in Neveu's setup, but allowed the players to use *randomized* stopping times: a strategy is a $[0,1]$-valued process that dictates the probability by which the player stops at every stage. They proved that the problem has a value, assuming only condition (a).

Extensive literature provides sufficient conditions for the existence of the value in continuous time [see, e.g., Bismut (1977), Alario-Nazaret, Lepeltier and Marchal (1982), Lepeltier and Maingueneau (1984), Touzi and Vieille (2002) and Laraki and Solan (2002)]. Some authors have studied the diffusion case, see, for example, Cvitanić and Karatzas (1996).

The nonzero–sum problem in discrete time when the payoffs have a special structure was studied, among others, by Mamer (1987), Morimoto (1986), Ohtsubo (1987, 1991), Nowak and Szajowski (1999) and Neumann, Ramsey and Szajowski (2002) and the references therein. In the nonzero–sum case, the processes $(R_{\{1\},n})$, $(R_{\{2\},n})$ and $(R_{\{1,2\},n})$ are $\mathbf{R}^2$-valued, and the expected payoff to player $i$, $i=1,2$, is

$$\gamma^i(\mu,\nu) = \mathbf{E}_{\mu,\nu}[\mathbf{1}_{\{\mu<\nu\}} R^i_{\{1\},\mu} + \mathbf{1}_{\{\mu>\nu\}} R^i_{\{2\},\nu} + \mathbf{1}_{\{\mu=\nu<\infty\}} R^i_{\{1,2\},\mu}].$$

The goal of each player is to maximize his or her own expected payoff. Given $\varepsilon > 0$, a pair of stopping times $(\mu,\nu)$ is an $\varepsilon$-*equilibrium* if for every pair of stopping times $(\mu',\nu')$,

$$\gamma^1(\mu,\nu) \geq \gamma^1(\mu',\nu) - \varepsilon \quad \text{and} \quad \gamma^2(\mu,\nu) \geq \gamma^2(\mu,\nu') - \varepsilon.$$

The above-mentioned authors provided various sufficient conditions under which $\varepsilon$-equilibria exist.

In the present paper, we study two-player nonzero–sum games in discrete time with randomized stopping times, and we prove the existence of an $\varepsilon$-equilibrium for every $\varepsilon > 0$, under merely an integrability condition. Our technique is based on a stochastic variation of Ramsey's theorem. Ramsey (1930) proved that for every coloring of a complete infinite graph by finitely many colors there is a complete infinite monochromatic subgraph. Our variation allows us to reduce the problem of the existence of an $\varepsilon$-equilibrium in a general stopping game to that of studying properties of $\varepsilon$-equilibria in a simple class of stochastic games with finite state space.

The paper is arranged as follows. In Section 2 we provide the model and the main result. A sketch of the proof appears in Section 3. In Section 4 we present a stochastic variation of Ramsey's theorem. In Section 5 we define the notion of games played on a finite tree and we study some of their properties. The proof of the main result in this section is relegated to Section 8. In Section 6 we show how to approximate a general filtration between two stopping times by a sequence of finite algebras. In Section 7 we construct an $\varepsilon$-equilibria. We end by discussing extensions to more than two players in Section 9.

Our proof uses tools both from the theory of stochastic processes (Sections 4 and 6) and from the theory of games (Sections 5 and 8). However, no prior knowledge of these fields is assumed.

**2. The model and the main result.** A *two-player nonzero–sum stopping game* is a 5-tuple $\Gamma = (\Omega, \mathcal{A}, p, \mathcal{F}, R)$ where:

- $(\Omega, \mathcal{A}, p)$ is a probability space;
- $\mathcal{F} = (\mathcal{F}_n)_{n \geq 0}$ is a filtration over $(\Omega, \mathcal{A}, p)$;
- $R = (R_n)_{n \geq 0}$ is an $\mathcal{F}$-adapted $\mathbf{R}^6$-valued process. The coordinates of $R_n$ are denoted by $R^i_{Q,n}$, $i = 1, 2$, $\phi \neq Q \subseteq \{1, 2\}$.

A (behavior) *strategy* for player 1 (resp. player 2) is a $[0,1]$-valued $\mathcal{F}$-adapted process $x = (x_n)_{n \geq 0}$ [resp. $y = (y_n)_{n \geq 0}$]. The interpretation is that $x_n$ (resp. $y_n$) is the probability by which player 1 (resp. player 2) stops at stage $n$ (provided the game is not stopped before that stage).

Let $\theta$ be the first stage, possibly infinite, in which at least one of the players stops and let $\phi \neq Q \subseteq \{1, 2\}$ be the set of players that stop at stage $\theta$ (provided $\theta < \infty$). The expected payoff under $(x, y)$ is given by

$$\gamma^i(x, y) = \mathbf{E}_{x,y}[R^i_{Q,\theta} \mathbf{1}_{\{\theta < \infty\}}], \quad (1)$$

where the expectation $\mathbf{E}_{x,y}$ is with respect to (w.r.t.) the distribution $\mathbf{P}_{x,y}$ over plays induced by $(x, y)$, and $\mathbf{1}$ is the indicator function.

DEFINITION 2.1. Let $\Gamma = (\Omega, \mathcal{A}, p, \mathcal{F}, R)$ be a two-player nonzero–sum stopping game and let $\varepsilon > 0$. A pair of strategies $(x^*, y^*)$ is an $\varepsilon$-*equilibrium* if $\gamma^1(x^*, y^*) \geq \gamma^1(x, y^*) - \varepsilon$ and $\gamma^2(x^*, y^*) \geq \gamma^2(x^*, y) - \varepsilon$ for every $x$ and $y$.

The main result of the paper is the following:

THEOREM 2.2. *Let $\Gamma = (\Omega, \mathcal{A}, p, \mathcal{F}, R)$ be a two-player nonzero–sum stopping game such that $\sup_{n \geq 0} \|R_n\|_\infty \in L^1(p)$. Then for every $\varepsilon > 0$, the game admits an $\varepsilon$-equilibrium.*

The definitions imply that for every two payoff processes $R$ and $\widetilde{R}$ such that $\mathbf{E}[\sup_{n \in \mathbf{N}} \|R_n - \widetilde{R}_n\|_\infty] < \varepsilon$, every $\varepsilon$-equilibrium in $(\Omega, \mathcal{A}, p, \mathcal{F}, \widetilde{R})$ is a $3\varepsilon$-equilibrium in $(\Omega, \mathcal{A}, p, \mathcal{F}, R)$. Hence we can assume without loss of generality (w.l.o.g.) that the payoff process $R$ is uniformly bounded and that its range is finite. Actually, we assume that for some $K \in \mathbf{N}$, $R_n \in \{0, \pm\frac{1}{K}, \pm\frac{2}{K}, \ldots, \pm\frac{K}{K}\}^6$ for every $n \in \mathbf{N}$.

**3. Sketch of the proof.** In the present section we provide the main ideas of the proof. Let $\Gamma$ be a stopping game. To simplify the presentation, assume that $\mathcal{F}_n$ is trivial for every $n$, so that the payoff process is deterministic. Recall that w.l.o.g. payoffs are uniformly bounded by 1.

Given $\varepsilon > 0$, fix a finite covering $M$ of the space of payoffs $[-1, 1]^2$ by sets with diameter smaller than $\varepsilon$. For every two nonnegative integers $k < l$, define the periodic game $G(k, l)$ to be the game that starts at stage $k$ and, if not stopped earlier, restarts at stage $l$. Formally, $G(k, l)$ is a stopping game in which the terminal payoff at stage $n$ is equal to the terminal payoff at stage $k + (n \bmod l - k)$ in $\Gamma$.

This periodic game is a simple stochastic game [see, e.g., Shapley (1953) or Flesch, Thvijsman and Vrieze (1996)] and is known to admit an $\varepsilon$-equilibrium in periodic strategies. Assign to each pair of nonnegative integers $k < l$ an element $m(k, l) \in M$ that contains the expected payoff that corresponds to a periodic $\varepsilon$-equilibria of the game $G(k, l)$.

Thus, we assigned to each $k < l$ a color $m(k, l) \in M$. A consequence of Ramsey's theorem is that there is an increasing sequence of integers $0 \leq k_1 < k_2 < \cdots$ such that $m(k_1, k_2) = m(k_n, k_{n+1})$ for every $n$.

Assume first that $k_1 = 0$. A naive candidate for a $3\varepsilon$-equilibrium suggests that between stages $k_n$ and $k_{n+1}$, the players follow a periodic $\varepsilon$-equilibrium in the game $G(k_n, k_{n+1})$ with corresponding payoff in the set $m(k_1, k_2)$.

For this strategy pair to indeed be a $3\varepsilon$-equilibrium, the properties of the $\varepsilon$-equilibria in periodic games must be studied. The complete solution of this case appears in Shmaya, Solan and Vieille (2003), who observed that in each periodic game $G(k, l)$ there exists a periodic $\varepsilon$-equilibrium that satisfies at least one of the following conditions: (i) Neither player ever stops. (ii) Both players receive nonnegative payoffs and termination occurs in each period with probability at least $\varepsilon^2$. (iii) If a player receives a negative payoff, then his or her opponent stops in each period with probability at least $\varepsilon^2$. The fact that at least one of these conditions holds is sufficient to prove that the concatenation described above is a $3\varepsilon$-equilibrium, with corresponding payoff in the convex hull of $m(k_1, k_2)$.

If $k_1 > 0$, choose an arbitrary $m \in m(k_1, k_2)$. Between stages 0 and $k_1$, the players follow an equilibrium in the $k_1$-stage game with terminal payoff $m$; that is, if no player ever stops before stage $k_1$, the payoff is $m$. From stage $k_1$ and on, the players follow the strategy described above. It is easy to verify that this strategy pair forms a $4\varepsilon$-equilibrium.

When the payoff process is general, few difficulties appear. First, a periodic game is defined now by two stopping times $\mu_1 < \mu_2$: $\mu_1$ indicates the initial stage and $\mu_2$ indicates when the game restarts. To analyze this periodic game, we have to reduce the problem to the case where the $\sigma$-algebras $\mathcal{F}_{\mu_1}, \mathcal{F}_{\mu_1+1}, \ldots, \mathcal{F}_{\mu_2}$ are finite. This is done in Section 7.

Second, we have to study properties of $\varepsilon$-equilibria in these periodic games, so that a proper concatenation of $\varepsilon$-equilibria in the different periodic games generates a $4\varepsilon$-equilibrium in the original game. This is done in Section 5.

Third, we have to generalize Ramsey's theorem to this stochastic setup. This is done in Section 4.

**4. A stochastic variation of Ramsey's theorem.** In the present section we provide a stochastic variation of Ramsey's theorem. Let $(\Omega, \mathcal{A}, p)$ be a probability space and let $\mathcal{F} = (\mathcal{F}_n)_{n \geq 0}$ be a filtration. For every set $A \subseteq \Omega$, $A^c = \Omega \setminus A$ is the complement of $A$. For every $A, B \in \mathcal{A}$, *A holds on B* if and only if $p(A^c \cap B) = 0$.

Ramsey (1930) proved that for every function that attaches a color $c(k,l) \in C$, where $C$ is a finite set, to every two nonnegative integers $k < l$, there is an increasing sequence of integers $k_0 < k_1 < \cdots$ such that $c(k_0, k_1) = c(k_i, k_j)$ for every $i < j$.

We are going to attach for every nonnegative integer $n$ and every stopping time $\tau$, an $\mathcal{F}_n$-measurable function $c_{n,\tau}$ that is defined over the set $\{\tau > n\}$, whose range is some finite set $C$. We also impose a consistency requirement: if $\tau_1 = \tau_2 > n$ on an $\mathcal{F}_n$-measurable set $F$, then $c_{n,\tau_1} = c_{n,\tau_2}$ on $F$. Under these conditions, we derive a weaker conclusion than that of Ramsey's theorem: for every $\varepsilon > 0$ there exists an increasing sequence of stopping times $\theta_0 < \theta_1 < \cdots$ such that $\mathbf{P}(c_{\theta_0, \theta_1} = c_{\theta_i, \theta_{i+1}} \, \forall \, i) > 1 - \varepsilon$.

We now formally present the result.

DEFINITION 4.1. An *NT function* is a function that assigns to every integer $n \geq 0$ and every bounded stopping time $\tau$ an $\mathcal{F}_n$-measurable random variable (r.v.) that is defined over the set $\{\tau > n\}$. We say that an NT function $f$ is *C-valued* for some set $C$ if the r.v. $f_{n,\tau}$ is $C$-valued for every $n \geq 0$ and every $\tau$.

DEFINITION 4.2. An NT function $f$ is $\mathcal{F}$-*consistent* if for every $n \geq 0$, every $\mathcal{F}_n$-measurable set $F$ and every two bounded stopping times $\tau_1, \tau_2$, we have

$$\tau_1 = \tau_2 > n \text{ on } F \qquad \text{implies} \qquad f_{n,\tau_1} = f_{n,\tau_2} \text{ on } F.$$

When $f$ is an NT function and $\sigma < \tau$ are two bounded stopping times, we denote $f_{\sigma,\tau}(\omega) = f_{\sigma(\omega),\tau}(\omega)$. Thus $f_{\sigma,\tau}$ is an $\mathcal{F}_\sigma$-measurable r.v.

The main result of this section is the following.

THEOREM 4.3. *For every finite set $C$ of colors, every $C$-valued $\mathcal{F}$-consistent NT function $c$ and every $\varepsilon > 0$, there exists a sequence of bounded stopping times $0 \leq \theta_0 < \theta_1 < \theta_2 < \cdots$ such that $p(c_{\theta_0,\theta_1} = c_{\theta_1,\theta_2} = c_{\theta_2,\theta_3} = \cdots) > 1 - \varepsilon$.*

COMMENT. The natural stochastic generalization of Ramsey's theorem requires the stronger condition $p(c_{\theta_0,\theta_1} = c_{\theta_i,\theta_j} \ \forall 0 \leq i < j) \geq 1 - \varepsilon$. We do not know whether this generalization is correct.

The following example shows that a sequence of stopping times $\theta_0 < \theta_1 < \theta_2 < \cdots$ such that $p(c_{\theta_0,\theta_1} = c_{\theta_1,\theta_2} = \cdots) = 1$ need not exist, even without the boundedness condition.

EXAMPLE 4.4. Let $X_n$ be a biased random walk on the integers, let $X_0 = 0$ and let $p(X_{n+1} = X_n + 1) = 1 - p(X_{n+1} = X_n - 1) = 3/4$. Let $\mathcal{F}_n = \sigma(X_0, X_1, \ldots, X_n)$. Let $R_0 = \Omega$ and, for every $n > 0$, let $R_n = \bigcup_{1 \leq k \leq n}\{X_k = -1\}$. For every finite (but not necessarily bounded) stopping time $\tau$ define $c_{n,\tau} = \text{Red}$ on $R_n \cap \{\tau > n\}$ and $c_{n,\tau} = \text{Blue}$ on $R_n^c \cap \{\tau > n\}$. Since $p(\bigcup_{n \geq 0} R_n) < 1$, whereas for every finite stopping time $\theta$ and every $B \in \mathcal{F}_\theta$, one has $p(\bigcup_{n \geq 0} R_n | B) > 0$, it follows that for every sequence $\theta_0 < \theta_1 < \cdots$ of finite stopping times $p(c_{\theta_0,\theta_1} = \text{Blue}) > 0$, whereas $p(c_{\theta_0,\theta_1} = c_{\theta_1,\theta_2} = \cdots = \text{Blue} | c_{\theta_0,\theta_1} = \text{Blue}) < 1$.

We start by proving a slightly stronger version of Theorem 4.3 for $|C| = 2$.

LEMMA 4.5. *Let $C = \{\text{Blue}, \text{Red}\}$ and let $c$ be a $C$-valued $\mathcal{F}$-consistent NT function. For every $\varepsilon > 0$ there exist $N \in \mathbf{N}$, two sets $\bar{R}, \bar{B} \in \mathcal{F}_N$ and a sequence $N \leq \tau_0 < \tau_1 < \tau_2 < \cdots$ of bounded stopping times, such that:*

(a) $\bar{R} = \bar{B}^c$;
(b) $p(c_{\tau_0,\tau_1} = c_{\tau_1,\tau_2} = \cdots = \text{Red} | \bar{R}) > 1 - \varepsilon$;
(c) $p(c_{\tau_k,\tau_l} = \text{Blue} \ \forall k, l | \bar{B}) > 1 - \varepsilon$.

PROOF. We claim first that for every $n \in \mathbf{N}$ we can find two sets $R_n, B_n \in \mathcal{F}_n$ and a bounded stopping time $\sigma_n$ such that:

1. We have $p(R_n \cup B_n) > 1 - 1/2^n$.
2. We have $\{\sigma_n > n\}$ on $R_n$ and $c_{n,\sigma_n} = \text{Red}$ on $R_n$.
3. For every bounded stopping time $\tau$, $c_{n,\tau} = \text{Blue}$ on $B_n \cap \{\tau > n\}$.

To see this, fix $n \in \mathbf{N}$. Call a set $F \in \mathcal{F}_n$ *red* if there exists a bounded stopping time $\sigma_F$ such that on $F$ both $\sigma_F > n$ and $c_{n,\sigma_F} = \text{Red}$. Since $c$ is $\mathcal{F}$-consistent, if $F, G \in \mathcal{F}_n$ are red, then so is $F \cup G$. Let $\alpha = \sup_F \{p(F), F \in \mathcal{F}_n \text{ is red}\}$. For every $k \geq 1$, let $F_k \in \mathcal{F}_n$ be a red set such that $p(F_k) > \alpha - \frac{1}{k}$. Let $F_* = \bigcup_{k \geq 1} F_k$. Observe that $F_* \in \mathcal{F}_n$ and $p(F_*) = \alpha$. Moreover, no subset of $F_*^c$ with positive probability is red. Let $R_n = F_{2^n}$, let $\sigma_n$ be a bounded stopping time such that on $R_n$, $\sigma_n > n$ and $c_{n,\sigma_n} = \text{Red}$, and let $B_n = F_*^c$. This concludes the proof of the claim.

Let $B = \{B_n \text{ i.o.}\}$, and set $R = B^c$. Since $R, B \in \bigvee_n \mathcal{F}_n$, $N \in \mathbf{N}$ and there are two sets $\bar{B}, \bar{R} \in \mathcal{F}_N$ such that (i) $\bar{R} = \bar{B}^c$, (ii) $p(B|\bar{B}) > 1 - \varepsilon$ and (iii) $p(R|\bar{R}) > 1 - \varepsilon$. On $R$, and therefore also on $R \cap \bar{R}$, both $B_n$ and $(B_n \cup R_n)^c$ occur only finitely many times. By sufficiently increasing $N$, we assume w.l.o.g. that $p(\bigcap_{n \geq N} R_n | R \cap \bar{R}) > 1 - \varepsilon$. In particular,

$$(2) \qquad p\left(\bigcap_{n \geq N} R_n \Big| \bar{R}\right) > 1 - 2\varepsilon.$$

Let $N = n_0 < n_1 < n_2 < \cdots$ be a sequence of integers such that, for every $k \geq 0$, $p(T_k | B \cap \bar{B}) > 1 - \varepsilon/2^{k+1}$, where $T_k = \bigcup_{n_k \leq n < n_{i+k}} B_n$. Then $p(\bigcap_{k \geq 0} T_k | B \cap \bar{B}) > 1 - \varepsilon$ and, therefore,

$$(3) \qquad p\left(\bigcap_{k \geq 0} T_k \Big| \bar{B}\right) > 1 - 2\varepsilon.$$

We now define the sequence $(\tau_k)_{k \geq 0}$ inductively, working separately on $\bar{R}$ and $\bar{B}$. Consider first the set $\bar{R}$. Define $\tau_0 = N$. Given $\tau_k$, define $\tau_{k+1} = \sum_{n \in \mathbf{N}} \sigma_n \mathbf{1}_{\{\tau_k = n\} \cap R_n \cap \bar{R}}$ on $\bar{R} \cap \bigcup_n (\{\tau_k = n\} \cap R_n)$. Since $\tau_k$ and $(\sigma_n)_{n \geq 0}$ are bounded, $\tau_{k+1}$ can be extended to a bounded stopping time on $\bar{R}$. By definition $c_{\tau_0, \tau_1} = c_{\tau_1, \tau_2} = \cdots = \text{Red}$ on $\bar{R} \cap (\bigcap_{n \geq N} R_n)$, and it follows from (2) that $p(c_{\tau_0, \tau_1} = c_{\tau_1, \tau_2} = \cdots = \text{Red} | \bar{R}) \geq 1 - 2\varepsilon$.

Consider now the set $\bar{B}$. Define $\tau_0 = N$. Define $\tau_{k+1}(w) = \min\{n_k \leq n < n_{k+1}, w \in B_n\}$ on $\bar{B} \cap T_k$ and $\tau_{k+1} = n_{k+1} - 1$ on $\bar{B} \setminus T_k$. By the definition of $\tau_k$, for every $k, l \in \mathbf{N}$, $c_{\tau_k, \tau_l} = \text{Blue}$ on $\bar{B} \cap (\bigcap_{k \geq 0} T_k)$ and it follows from (3) that $p(c_{\tau_k, \tau_l} = \text{Blue } \forall k, l | \bar{B}) > 1 - 2\varepsilon$. $\square$

PROOF OF THEOREM 4.3. We prove the theorem by induction on $|C|$. The case $|C| = 2$ follows from Lemma 4.5. Assume we have already proven the theorem for $|C| = r$ and assume $|C| = r + 1$. Let Red be a color in $C$.

By considering all colors other than Red as a single fictitious color and by applying Lemma 4.5, there exist $N \in \mathbf{N}$, two sets $\bar{R}, \bar{B} \in \mathcal{F}_N$ and a sequence of stopping times $N \leq \tau_0 < \tau_1 < \cdots$ such that (i) $\bar{R} = \bar{B}^c$, (ii) $p(c_{\tau_0, \tau_1} = c_{\tau_1, \tau_2} = \cdots = \text{Red} | \bar{R}) > 1 - \varepsilon/2$ and (iii) $p(c_{\tau_k, \tau_l} \neq \text{Red } \forall k, l | \bar{B}) > 1 - \varepsilon/2$. We define $\theta_i$ separately on $\bar{R}$ and $\bar{B}$. On $\bar{R}$, we let $\theta_i = \tau_i$.

We now restrict ourselves to the space $(\bar{B}, \mathcal{A}_{\bar{B}}, p_{\bar{B}})$ with the filtration $\mathcal{G}_n = \mathcal{F}_{\tau_n} \cap \bar{B}$. Let $\tilde{c}$ be the $C$-valued NT function over $\mathcal{G}$ defined by $\tilde{c}_{n,\beta} = c_{\tau_n, \tau_\beta}$ for every stopping time $\beta$ of $\mathcal{G}$, where $\tau_\beta = \sum_n \tau_n \mathbf{1}_{\{\beta = n\}}$ is a stopping time of $\mathcal{F}$. Let $c'$ be the coloring that is obtained from $\tilde{c}$ by swapping the color Red with another color in $C$, say Green:

$$c'_{n,\beta} = \begin{cases} \tilde{c}_{n,\beta}, & \text{if } \tilde{c}_{n,\beta} \neq \text{Red}, \\ \text{Green}, & \text{if } \tilde{c}_{n,\beta} = \text{Red}. \end{cases}$$

Since $c'$ is a $C \setminus \{\text{Red}\}$-valued $\mathcal{G}$-consistent NT function, we can apply the induction hypothesis and obtain a sequence of stopping times $0 \leq \beta_0 < \beta_1 < \beta_2 < \cdots$ of $\mathcal{G}$ such that

(4) $$p(c'_{\beta_0,\beta_1} = c'_{\beta_1,\beta_2} = \cdots | \bar{B}) > 1 - \varepsilon/2.$$

By (4) and (iii) it follows that $p(\tilde{c}_{\beta_0,\beta_1} = \tilde{c}_{\beta_1,\beta_2} = \cdots | \bar{B}) > 1 - \varepsilon$. We define $\theta_i = \tau_{\beta_i}$ on $\bar{B}$. Thus

(5) $$p(c_{\theta_0,\theta_1} = c_{\theta_1,\theta_2} = \cdots | \bar{B}) > 1 - \varepsilon.$$

Combining (ii) and (5) we get $p(c_{\theta_0,\theta_1} = c_{\theta_1,\theta_2} = \cdots) > 1 - \varepsilon$, as desired. □

**5. Stopping games on finite trees.** An important building block in our analysis is stopping games that are played on a finite tree. In the present section we define these games and study some of their properties.

5.1. *The model.*

DEFINITION 5.1. A *stopping game on a finite tree* (or simply a *game on a tree*) is a tuple $T = (S, S_1, r, (C_s, p_s, R_s)_{s \in S \setminus S_1})$, where

- $(S, S_1, r, (C_s)_{s \in S \setminus S_1})$ is a tree, $S$ is a nonempty finite set of *nodes*, $S_1 \subseteq S$ is a nonempty set of *leaves*, $r \in S$ is the root and, for each $s \in S \setminus S_1$, $C_s \subseteq S \setminus \{r\}$ is the nonempty set of *children* of $s$ [we denote by $S_0 = S \setminus S_1$ the set of nodes which are not leaves; for every $s \in S$, depth($s$) is the depth of $s$—the length of the path that connects the root to $s$] and for every $s \in S_0$,

- $p_s$ is a probability distribution over $C_s$;
- $R_s \in \mathbf{R}^6$ is the *payoff* at $s$ [the coordinates of $R_s$ are denoted $(R^i_{Q,s})_{i=1,2,\ \phi \neq Q \subseteq \{1,2\}}$].

A stopping game on a finite tree starts at the root and is played in stages. Given the current node $s \in S_0$, and the sequence of nodes already visited, both players decide, simultaneously and independently, whether to stop or to continue. Let $Q$ be the set of players that decide to stop. If $Q \neq \phi$, the play terminates and the terminal payoff to each player $i$ is $R^i_{Q,s}$. If $Q = \phi$, a new node $s'$ in $C_s$ is chosen according to $p_s$. The process now repeats itself, with $s'$ being the current node. If $s' \in S_1$, the new current node is the root $r$. Thus, players cannot stop at leaves.

The game on the tree is essentially played in rounds. The round starts at the root and ends once it reaches a leaf.

Consider the first round of the game. Let $t$ be the stopping stage. If no termination occurs in the first round, $t = \infty$. If $t < \infty$, let $s$ be the node (of depth $t$) in which termination occurred and let $Q$ be the set of players that stop at stage $t$. The r.v. $r^i = R^i_{Q,s} \mathbf{1}_{\{t<\infty\}}$ is the payoff to player $i$ in the first round.

A *stationary strategy* of player 1 (resp. player 2) is a function $x: S_0 \to [0,1]$ (resp. $y: S_0 \to [0,1]$); $x(s)$ is the probability that player 1 stops at $s$. Denote by $\mathbf{P}_{x,y}$ the distribution over plays induced by $(x,y)$, and denote by $\mathbf{E}_{x,y}$ the corresponding expectation operator.

For every pair of stationary strategies $(x,y)$ we denote by $\pi(x,y) = \mathbf{P}_{x,y}(t < \infty)$ the probability that under $(x,y)$ the game terminates in the first round of the game; that is, the probability that the root is visited only once along the play. We denote by $\rho^i(x,y) = \mathbf{E}_{x,y}[r^i]$, $i=1,2$, the expected payoff of player $i$ in a single round. Finally, we set $\gamma^i(x,y) = \rho^i(x,y)/\pi(x,y)$ (by convention, $\frac{0}{0} = 0$). This is the expected payoff under $(x,y)$. In particular,

$$\pi(x,y) \times \gamma^i(x,y) = \rho^i(x,y). \tag{6}$$

When we want to emphasize the dependency of these variables on the game $T$, we write $\pi(x,y;T)$, $\rho^i(x,y;T)$ and $\gamma^i(x,y;T)$.

Observe that for every pair of stationary strategies $(x,y)$,

$$\pi(x,0) + \pi(0,y) \geq \pi(x,y), \tag{7}$$

where 0 is the strategy that never stops; that is, $0(s) = 0$ for every $s$.

DEFINITION 5.2. A pair of stationary strategies $(x,y)$ is an $\varepsilon$-*equilibrium* of the game on a tree $T$ if, for each pair of strategies $(x',y')$, $\gamma^1(x',y) \leq \gamma^1(x,y) + \varepsilon$ and $\gamma^2(x,y') \leq \gamma^2(x,y) + \varepsilon$.

COMMENT. A stopping game on a finite tree $T$ is equivalent to a recursive absorbing game, where each round of the game $T$ corresponds to a single stage of the recursive absorbing game. A recursive absorbing game is a stochastic game with a single nonabsorbing state in which the payoff in nonabsorbing states is 0. Flesch, Thuijsman and Vrieze (1996) proved that every recursive absorbing game admits an $\varepsilon$-equilibrium in stationary strategies. This result also follows from the analysis of Vrieze and Thuijsman (1989). However, there is no bound on the per-round probability of termination under this $\varepsilon$-equilibrium and we need to bound this quantity.

5.2. *Main results concerning games on trees.* Throughout this section we fix $R^1, R^2 \in \mathbf{R}$, such that at least one of them is positive, and play a game on a tree whose payoffs $(R_s)_{s \in S_0}$ satisfy the following conditions for every $i = 1, 2$, every $\varnothing \subset Q \subseteq \{1, 2\}$ and every node $s \in S_0$:

B1. For some $K \in \mathbf{N}$, $R^i_{Q,s} \in \{0, \pm\frac{1}{K}, \ldots, \pm\frac{K}{K}\}$.
B2. We have $R^i_{\{i\},s} \leq R^i$.
B3. Whenever $R^1_{\{1\},s} = R^1$, $R^2_{\{1\},s} < R^2$.
B4. Whenever $R^1_{\{2\},s} = R^2$, $R^1_{\{2\},s} < R^1$.

We have already seen that condition B1 can be assumed w.l.o.g. We will later set $R^i$ to be an upper bound of $(R^i_{\{i\},n})$, so that condition B2 can be assumed. The results we prove in this section are not trivial only when $R^1_{\{1\},s} = R^1$ and $R^2_{\{2\},s'} = R^2$ for some $s, s' \in S_0$. As we see later, when condition B3 or B4 does not hold, a simple $\varepsilon$-equilibrium exists.

Assuming no player ever stops, the collection $(p_s)_{s \in S_0}$ of probability distributions at the nodes induces a probability distribution over the set $S_1$ of leaves or, equivalently, over the set of branches that connect the root to the leaves. For each set $D \subseteq S_0$, we denote by $p_D$ the probability that the chosen branch passes through $D$. For each $s \in S$, we denote by $F_s$ the event that the chosen branch passes through $s$.

We first bound the probability of termination in a single round when the $\varepsilon$-equilibrium payoff is low for at least one player.

LEMMA 5.3. *Let $\varepsilon > 0$ and let $(x, y)$ be a stationary $\frac{\varepsilon}{2}$-equilibrium in $T$ such that $\gamma^1(x, y) \leq R^1 - \varepsilon$. Then $\pi(0, y) \geq \frac{\varepsilon}{6} \cdot \mu_1$, where $\mu_1 = \mu_1(T) = p(\bigcup \{F_s, R^1_{\{1\},s} = R^1\})$ is the probability that if both players never stop, the game visits a node $s$ with $R^1_{\{1\},s} = R^1$ in the first round.*

An analogous statement holds for player 2.

PROOF OF LEMMA 5.3. Consider the following strategy $z$ of player 1:

$$z_s = \begin{cases} 1, & \text{if } R^1_{\{1\},s} = R^1, \\ 0, & \text{otherwise.} \end{cases}$$

Denote by $t_1$ and $t_2$ the stopping stages of the two players in the first round. By the definition of $z$, and since payoffs are bounded by 1,

$$\rho^1(z, y) = \mathbf{P}_{z,y}(t_1 < \min\{\infty, t_2\}) \times R^1 + \mathbf{E}_{z,y}[\mathbf{1}_{t_2 \leq t_1} r^1]$$
(8)
$$\geq \mathbf{P}_{z,y}(t_1 < \infty) \times R^1 - 2\mathbf{P}_{z,y}(t_2 < \infty)$$
$$= \mu_1 R^1 - 2\pi(0, y).$$

Since $(x,y)$ is an $\frac{\varepsilon}{2}$-equilibrium and since, by (7), $\pi(z,y) \leq \pi(0,y) + \pi(z,0) = \pi(0,y) + \mu_1$,

$$\begin{aligned}\rho^1(z,y) &= \gamma^1(z,y) \times \pi(z,y) \\ &\leq \left(\gamma^1(x,y) + \frac{\varepsilon}{2}\right) \times (\pi(0,y) + \mu_1) \\ &\leq \left(R^1 - \frac{\varepsilon}{2}\right) \times (\pi(0,y) + \mu_1).\end{aligned} \quad (9)$$

Equations (8) and (9) imply that $\pi(0,y) \geq \frac{\varepsilon}{6} \times \mu_1$. $\square$

DEFINITION 5.4. Let $T = (S, S_1, r, (C_s, p_s, R_s)_{s \in S_0})$ and $T' = (S', S'_1, r', (C'_s, p'_s, R'_s)_{s \in S'_0})$ be two games on trees. We say $T'$ is a *subgame* of $T$ if (i) $S' \subseteq S$, (ii) $r' = r$ and (iii) for every $s \in S'_0$, $C'_s = C_s$, $p'_s = p_s$ and $R'_s = R_s$.

In words, $T'$ is a subgame of $T$ if we remove all the descendants (in the strict sense) of several nodes from the tree $(S, S_1, r, (C_s)_{s \in S_0})$ and keep all other parameters fixed. Observe that this notion is different from the standard definition of a subgame in game theory.

Let $T = (S, S_1, r, (C_s, p_s, R_s)_{s \in S_0})$ be a game on a tree. For each subset $D \subseteq S_0$, we denote by $T_D$ the subgame of $T$ generated by trimming $T$ from $D$ downward. Thus, all strict descendants of nodes in $D$ are removed.

For every subgame $T'$ of $T$ and every subgame $T''$ of $T'$, let $p_{T'',T'} = p_{S''_1 \setminus S'_1}$ be the probability that the chosen branch in $T$ passes through a leaf of $T''$ strictly before it passes through a leaf of $T'$ (here, $S'_1$ and $S''_1$ are the sets of leaves of $T'$ and $T''$, resp.).

The next proposition analyzes $\varepsilon$-equilibria that yield a high payoff to both players. Since its proof is involved and independent of the rest of the paper, it is deferred to Section 8.

PROPOSITION 5.5. *Let $\varepsilon \in (0, 1/(36K^2))$ and, for $i \in \{1,2\}$, let $a_i \geq R^i - \varepsilon$. There exist a set $D \subseteq S_0$ of nodes and a strategy pair $(x,y)$ in $T$ such that:*

1. *No subgame of $T_D$ has an $\varepsilon$-equilibrium with corresponding payoffs in $[a_1, a_1 + \varepsilon] \times [a_2, a_2 + \varepsilon]$.*
2. *Either (a) $D = \phi$ (so that $T_D = T$) or (b) $(x,y)$ is a $9\varepsilon$-equilibrium in $T$, $a_i - \varepsilon \leq \gamma^i(x,y)$ and $\pi(x,y) \geq \varepsilon^2 \times p_D$.*

COMMENT. Actually we prove that in case 2(b), for every pair $(x', y')$ of strategies, $\gamma^1(x', y) \leq a_1 + 8\varepsilon$ and $\gamma^2(x, y') \leq a_2 + 8\varepsilon$.

5.3. *Coloring a finite tree.* In the present section we provide an algorithm that for every finite tree $T$ attaches a color $c(T)$ and several numbers $(\lambda_j(T))_j$ in the unit interval.

A rectangle $[a_1, a_1 + \varepsilon] \times [a_2, a_2 + \varepsilon]$ is *bad* if $R^1 - \varepsilon \leq a_1$ and $R^2 - \varepsilon \leq a_2$. It is *good* if $a_1 + \varepsilon \leq R^1 - \varepsilon$ or $a_2 + \varepsilon \leq R^2 - \varepsilon$.

Let $M$ be a finite covering of $[-1,1]^2$ with (not necessarily disjoint) rectangles $[a_1, a_1 + \varepsilon] \times [a_2, a_2 + \varepsilon]$, all of which are either good or bad. Thus, for every $u \in [-1,1]^2$ there is a rectangle $m \in M$ such that $u \in m$. We denote by $H = \{h_1, h_2, \ldots, h_J\}$ the set of bad rectangles in $M$ and denote by $G = \{g_1, g_2, \ldots, g_V\}$ the set of good rectangles in $M$.

Set $C = G \cup \{\varnothing\}$. This set is composed of the set $G$ of good rectangles together with another symbol $\varnothing$. For every game on a tree $T$ consider the following procedure which attaches an element $c \in C$ to $T$:

- Set $T^{(0)} = T$.
- For $1 \leq j \leq J$ apply Proposition 5.5 to $T^{(j-1)}$ and the bad rectangle $h_j = [a_{j,1}, a_{j,1} + \varepsilon] \times [a_{j,2}, a_{j,2} + \varepsilon]$ to obtain a subgame $T^{(j)}$ of $T^{(j-1)}$ and strategies $(x_T^{(j)}, y_T^{(j)})$ in $T^{(j)}$ such that:

  1. No subgame of $T^{(j)}$ has an $\varepsilon$-equilibrium with corresponding payoffs in $h_j$.
  2. Either $T^{(j)} = T^{(j-1)}$ or the following three conditions hold:
     (a) For $i \in \{1, 2\}$, $a_{j,i} - \varepsilon \leq \gamma^i(x_T^{(j)}, y_T^{(j)})$
     (b) For every pair $(x', y')$, $\gamma^1(x', y_T^{(j)}) \leq a_{j,1} + 8\varepsilon$ and $\gamma^2(x_T^{(j)}, y') \leq a_{j,2} + 8\varepsilon$.
     (c) We have $\pi(x_T^{(j)}, y_T^{(j)}) \geq \varepsilon^2 \times p_{T^{(j)}, T^{(j-1)}}$, where $p_{T^{(j)}, T^{(j-1)}}$ is the probability that a randomly chosen branch passes through a leaf of $T^{(j)}$ which is not a leaf of $T^{(j-1)}$ (see Section 5.1).

- If $T^{(J)}$ is trivial (i.e., the only node is the root), set $c(T) = \varnothing$; otherwise choose a stationary $\frac{\varepsilon}{2}$-equilibrium $(x^{(0)}, y^{(0)})$ of $T^{(J)}$. Since $T^{(J)}$ is a subgame of each $T^{(j)}$, $j = 1, \ldots, J$, since no subgame of $T^{(j)}$ has an $\varepsilon$-equilibrium with corresponding payoffs in $h_j$ and since every $\frac{\varepsilon}{2}$-equilibrium is also an $\varepsilon$-equilibrium, the corresponding $\frac{\varepsilon}{2}$-equilibrium payoff lies in some good rectangle $g \in G$. Set $c(T) = g$.

The strategies $(x_T^{(j)}, y_T^{(j)})$, as given by Proposition 5.5, are strategies in $T^{(j-1)}$. We consider them as strategies in $T$ by letting them continue from the leaves of $T^{(j-1)}$ downward.

We also define, for every $j \in J$,
$$\lambda_j(T) = p_{T^{(j)}, T^{(j-1)}},$$
so that $\pi(x_T^{(j)}, y_T^{(j)}) \geq \varepsilon^2 \times \lambda_j(T)$.

**6. Representative approximations.** Theorem 4.3 enables us to reduce the analysis to finite-stage games: games that start at some stage $n$ and terminate at stage $\tau$, where $\tau$ is some bounded stopping time. Since the state space $\Omega$ is arbitrary, while our game-theoretic tools allow us to analyze only games that are defined over a finite state space, we need to approximate each $\mathcal{F}_n$ by a finite $\sigma$-algebra.

Roughly, our goal here is to define a consistent NT function $T = (T_{n,\tau})$ that assigns for every $n \geq 0$, every bounded stopping time $\tau$ and every $\omega \in \Omega$, a game on a tree $T_{n,\tau}(\omega)$ that approximates the finite-stage game that is played between stages $n$ and $\tau$ in a desirable way.

This is done in two stages. First, we define, for every $k \geq n$, a finite partition $\mathcal{G}_k$ of $\Omega$ such that $(\mathcal{G}_k)_{k \geq n}$ contains all the information relevant to the players between stages $n$ and $\tau$. Second, given the sequence of finite partitions, we define the games on a tree $(T_{n,\tau}(\omega))_{\omega \in \Omega}$.

Throughout this section we fix a stopping game $\Gamma$ and $\varepsilon > 0$. Denote $\delta_n = \varepsilon^2/2^{n+2}$ for each $n \geq 0$. Set $\Delta_n = \sum_{k \geq n} \delta_k = \varepsilon^2/2^{n+1}$, so that $\sum_{n \geq 0} \Delta_n = \varepsilon^2$.

6.1. *Partial games.* In the present section we consider the partial game that is played between stages $\tau_1$ and $\tau_2$, where $\tau_1 \leq \tau_2$ are two bounded stopping times, and we define and study the notion of approximating games on a tree.

DEFINITION 6.1. Let $\tau_1 \leq \tau_2$ be two bounded stopping times. An $(\mathcal{F}, \tau_1, \tau_2)$-*strategy* is a sequence $x = (x_k)$ of random variables such that for every $k \geq 0$, (i) $x_k : \{\tau_1 \leq k < \tau_2\} \to [0,1]$ and (ii) $x_k$ is $\mathcal{F}_k$-measurable.

Thus, an $(\mathcal{F}, \tau_1, \tau_2)$-strategy prescribes to the player what to play between stages $\tau_1$ and $\tau_2$ (excluded). If $0 = \tau_0 < \tau_1 < \cdots$ is an increasing sequence of bounded stopping times and if, for each $l \geq 0$, $x^l$ is an $(\mathcal{F}, \tau_l, \tau_{l+1})$-strategy, we can naturally define a strategy $x$ in the stopping game simply by concatenating the strategies $(x^l)_{l \geq 0}$. Similarly, any strategy $x = (x_k)_{k \geq 0}$ in the stopping game naturally defines an $(\mathcal{F}, \tau_1, \tau_2)$-strategy for every pair of bounded stopping times $\tau_1 \leq \tau_2$ by considering the proper restriction of $x$.

When $x$ and $y$ are two $(\mathcal{F}, \tau_1, \tau_2)$-strategies, we define by $\pi(x, y; \mathcal{F}, \tau_1, \tau_2) = \mathbf{P}_{x,y}[\mathbf{1}_{\{\tau_1 \leq \theta < \tau_2\}}|\mathcal{F}_{\tau_1}]$ the conditional probability under $(x, y)$ that the game that starts at stage $\tau_1$ ends before stage $\tau_2$, and define by

$$\rho(x, y; \mathcal{F}, \tau_1, \tau_2) = \mathbf{E}_{x,y}[\mathbf{1}_{\{\tau_1 \leq \theta < \tau_2\}} R_{Q,\theta}|\mathcal{F}_{\tau_1}]$$

the corresponding expected payoff. We define

$$\gamma(x, y; \mathcal{F}, \tau_1, \tau_2) = \frac{\rho(x, y; \mathcal{F}, \tau_1, \tau_2)}{\pi(x, y; \mathcal{F}, \tau_1, \tau_2)}.$$

These are $\mathcal{F}_{\tau_1}$-measurable r.v.'s.

DEFINITION 6.2. Let $\tau_1 \leq \tau_2$ be two bounded stopping times. A $\delta$ *approximation* of $\Gamma$ between $\tau_1$ and $\tau_2$ is a pair $((\mathcal{G}_k), (q_{G,k}))$ such that for every $k \geq 0$:

1. $\mathcal{G}_k$ is an $\mathcal{F}_k$-measurable finite partition of $\{\tau_1 \leq k \leq \tau_2\}$;
2. $R_k$ is $\mathcal{G}_k$-measurable;
3. $\tau_1$ and $\tau_2$ are measurable w.r.t. $\mathcal{G}$; that is, for every $k \geq 0$, $\{\tau_1 = k\}$ and $\{\tau_2 = k\}$ are unions of atoms in $\mathcal{G}_k$;
4. any atom $G$ of $\mathcal{G}_k$ such that $k < \tau_2$ on $G$ is a union of some atoms in $\mathcal{G}_{k+1}$;
5. for every atom $G$ of $\mathcal{G}_k$, $q_{G,k}$ is a probability distribution over the atoms of $\mathcal{G}_{k+1}$ that are contained in $G$;
6. $\sum_{G' \in \mathcal{G}_{k+1}} |\mathbf{P}(G'|\mathcal{F}_k)(\omega) - q_{G,k}(G')| < \delta_k$ for every atom $G$ of $\mathcal{G}_k$ and almost every $\omega \in G$.

We identify $\mathcal{G}_k$ with the finite $\sigma$-algebra generated by $\mathcal{G}_k$ and we denote $\mathcal{G} = (\mathcal{G}_k)_{k \geq 0}$.

With every $\delta$ approximation $(\mathcal{G}, (q_{G,k}))$ of $\Gamma$ between $\tau_1$ and $\tau_2$ and every atom $G$ of $\mathcal{G}_{\tau_1}$ we can attach a game on a tree.

- The root is $G$.
- The nodes are all nonempty atoms $F$ of $(\mathcal{G}_k)$ such that (a) $F \subseteq G$ and (b) if $F \in \mathcal{G}_k$, then $\tau \geq k$ on $F$.
- The leaves are all atoms $F \in \bigcup_{k \geq n} \mathcal{G}_k$ where there is equality in (b).
- The payoff is given by $(R_k)_{\tau_1 \leq k \leq \tau_2}$.
- The children of each atom $F$ in $\mathcal{G}_k$ are all atoms $F'$ in $\mathcal{G}_{k+1}$ which are subsets of $F$.
- The transition from any node $F$ in $\mathcal{G}_k$ is given by $q_{F,k}$.

We denote this game on a tree by $T(\tau_1, \tau_2; \mathcal{G}, (q_{G,k}), G)$.

Suppose that for every atom $G$ of $\mathcal{G}_{\tau_1}$, $x_G$ is a strategy in the game $T(\tau_1, \tau_2; \mathcal{G}, (q_{G,k}), G)$. The collection of those strategies naturally defines a $(\mathcal{F}, \tau_1, \tau_2)$-strategy.

Similarly, if $x = (x_k)$ is an $(\mathcal{F}, \tau_1, \tau_2)$-strategy such that for every $k$, $x_k$ is $\mathcal{G}_k$-measurable, then $x$ naturally defines a strategy in $T(\tau_1, \tau_2; \mathcal{G}, (q_{G,k}), G)$ for every atom $G$ of $\mathcal{G}_{\tau_1}$.

Let $(x, y)$ be a pair of $(\mathcal{F}, \tau_1, \tau_2)$-strategies such that for every $k$, $x_k$ and $y_k$ are $\mathcal{G}_k$-measurable. We denote by $\pi(x, y; \tau_1, \tau_2, \mathcal{G}, (q_{G,k}))$ the probability of termination in one round of the game on a tree $T(\tau_1, \tau_2; \mathcal{G}, (q_{G,k}), G)$ under the strategies $(x, y)$, by $\rho(x, y; \tau_1, \tau_2, \mathcal{G}, (q_{G,k}))$ the expected payoff in one round, and by

$$\gamma(x, y; \tau_1, \tau_2, \mathcal{G}, (q_{G,k})) = \frac{\rho(x, y; \tau_1, \tau_2, \mathcal{G}, (q_{G,k}))}{\pi(x, y; \tau_1, \tau_2, \mathcal{G}, (q_{G,k}))}$$

the expected payoff in the game. These three functions are $\mathcal{G}_{\tau_1}$-measurable.

The following lemma provides estimates for the difference between the expected payoff and the expected probability of termination, when the filtration is changed. Its proof is omitted.

LEMMA 6.3. *Let $\tau_1 \leq \tau_2$ be two bounded stopping times and let $(\mathcal{G}, (q_{G,k}))$ be a $\delta$ approximation of $\Gamma$ between $\tau_1$ and $\tau_2$. Let $x$ and $y$ be a pair of $(\mathcal{G}, \tau_1, \tau_2)$-strategies. Then:*

1. $|\rho^i(x,y;\mathcal{F},\tau_1,\tau_2) - \rho^i(x,y;\tau_1,\tau_2,\mathcal{G},(q_{G,k}))| < \Delta_{\tau_1}$ *for each $i = 1, 2$.*
2. $|\pi(x,y;\mathcal{F},\tau_1,\tau_2) - \pi(x,y;\tau_1,\tau_2,\mathcal{G},(q_{G,k}))| < \Delta_{\tau_1}$.

The following lemma states that if $(\mathcal{G}, (q_{G,k}))$ is a $\delta$ approximation of $\Gamma$ between $\tau_1$ and $\tau_2$, and if the opponent plays a $(\mathcal{G}, \tau_1, \tau_2)$-strategy, then the player does not lose much by considering only $(\mathcal{G}, \tau_1, \tau_2)$-strategies [rather than $(\mathcal{F}, \tau_1, \tau_2)$-strategies].

LEMMA 6.4. *Let $\tau_1 \leq \tau_2$ be two bounded stopping times and let $(\mathcal{G}, (q_{G,k}))$ be a $\delta$ approximation of $\Gamma$ between $\tau_1$ and $\tau_2$. Let $x$ be a $(\mathcal{G}, \tau_1, \tau_2)$-strategy for player 1 and set $\gamma = \operatorname{ess\,sup}\{\gamma^2(x,y;\mathcal{G},\tau_1,\tau_2), y \text{ is a } (\mathcal{G}, \tau_1, \tau_2)\text{-strategy}\}$. Then, for every $(\mathcal{F}, \tau_1, \tau_2)$-strategy $y$,*

$$\rho^2(x,y;\mathcal{F},\tau_1,\tau_2) \leq \gamma \times \pi(x,y;\mathcal{F},\tau_1,\tau_2) + \Delta_{\tau_1} \qquad a.e.$$

PROOF. For $\mathcal{H} = \mathcal{F}, \mathcal{G}$, let

$$\alpha(\mathcal{H}) = \operatorname{ess\,sup}\{\rho^2(x,y;\mathcal{H},\tau_1,\tau_2) \\ + \gamma \times (1 - \pi(x,y;\mathcal{H},\tau_1,\tau_2)), y \text{ is a } (\mathcal{H},\tau_1,\tau_2)\text{-strategy}\}.$$

When player 2 is restricted to use $(\mathcal{H}, \tau_1, \tau_2)$-strategies, $\alpha(\mathcal{H})$ is his or her best possible payoff in the game that starts at stage $\tau_1$ and, if no player stops before stage $\tau_2$, terminates with payoff $\gamma$. From the definition of $\gamma$ it follows that $\alpha(\mathcal{G}) \leq \gamma$. Define recursively

$$(10) \quad \alpha(k, \mathcal{H}) = \begin{cases} \gamma, & k \geq \tau_2 \\ \max\{x_k \times R^2_{\{1\},k} + (1-x_k) \times \mathbf{E}(\alpha(k+1,\mathcal{H})|\mathcal{H}_k), \\ \qquad x_k \times R^2_{\{1,2\},k} + (1-x_k) \times R^2_{\{2\},k}\}, \\ & k < \tau_2. \end{cases}$$

The first term of the maximization in the second line corresponds to player 2 not stopping at stage $k$, while the second term corresponds to player 2 stopping at that stage. Plainly $\alpha(\mathcal{H}) = \alpha(\tau_1, \mathcal{H})$.

Since $x_k$ and $R^2_k$ are $\mathcal{G}_k$-measurable and since $(\mathcal{G}, (q_{G,k}))$ is a $\delta$ approximation of $\Gamma$ between $\tau_1$ and $\tau_2$, it follows by induction that $\alpha(k, \mathcal{F}) \leq \alpha(k, \mathcal{G}) +$ $\sum_{j=k}^{\infty} \delta_j$ for every $k \geq 0$. In particular, $\alpha(\mathcal{F}) = \alpha(\tau_1, \mathcal{F}) \leq \alpha(\tau_1, \mathcal{G}) + \Delta_{\tau_1} \leq \gamma + \Delta_{\tau_1}$. It follows that for every $(\mathcal{F}, \tau_1, \tau_2)$-strategy $y$, $\rho^2(x, y; \mathcal{F}, \tau_1, \tau_2) + \gamma \times (1 - \pi(x, y; \mathcal{F}, \tau_1, \tau_2)) \leq \gamma + \Delta_{\tau_1}$, which implies $\rho^2(x, y; \mathcal{F}, \tau_1, \tau_2) \leq \gamma \times \pi(x, y; \mathcal{F}, \tau_1, \tau_2) + \Delta_{\tau_1}$. $\square$

6.2. *Approximating the filtration $\mathcal{F}$.* Our main result in this section is the following:

THEOREM 6.5. *Let $\Gamma$ be a two-player nonzero–sum stopping game. There is a consistent NT function that assigns to every $n \geq 0$ and every bounded stopping time $\tau$ a $\delta$ approximation of $\Gamma$ between $n$ and $\tau$.*

The proof of this result, though quite technical, is intuitive. We start at stage $\tau$ and we proceed backward until we reach stage $n$. For every stage $n \leq k \leq \tau$, we define a finite partition of $\Omega$. Roughly, the partition $\mathcal{G}_k$ at stage $k$ is defined in such a way that (i) payoffs at stage $k$ are measurable w.r.t. $\mathcal{G}_k$ and (ii) the conditional probability to reach any atom $G \in \mathcal{G}_{k+1}$ does not vary by much over each atom of $\mathcal{G}_k$. The construction is slightly complicated since it should be consistent.

PROOF OF THEOREM 6.5. For every $i \geq 0$ and every $n \in \mathbf{N}$, choose once and for all a partition $\mathcal{B}_i^n$ of the $(n-1)$-dimensional simplex $\{r \in \mathbf{R}^n : \sum_{j=1}^n r_j = 1, r_j \geq 0 \ \forall j\}$ such that the diameter of each element in $\mathcal{B}_i^n$ is less than $\delta_i$ in the norm $\|\cdot\|_1$. For each $B \in \mathcal{B}_i^n$ choose an element $q_B \in B$.

Fix a total order on the collection of subsets of $\Omega$. This enables us to identify, for every partition $\mathcal{G}$ of $\Omega$ to $n$ atoms, the space $\mathbf{R}^n$ with the space $\mathbf{R}^{\mathcal{G}}$, simply by identifying the $i$th coordinate of $\mathbf{R}^n$ with the $i$th atom of $\mathcal{G}$ according to the total order.

Fix a nonnegative integer $n \geq 0$ and a bounded stopping time $\tau$. Let $M$ be a bound on $\tau$. We first define a sequence of finite partitions $(\tilde{\mathcal{G}}_k)$ recursively from $M$ down to $n$. We say that $\omega \approx_M \omega'$ if and only if (i) $\tau(\omega) = \tau(\omega') = M$ and (ii) $R_M(\omega) = R_M(\omega')$. We let $\tilde{\mathcal{G}}_M$ be the partition of $\Omega$ induced by the $\approx_M$ relation. Then $\tilde{\mathcal{G}}_M \subseteq \mathcal{F}_M$.

Suppose we have already defined a finite partition $\tilde{\mathcal{G}}_{m+1} \subseteq \mathcal{F}_{m+1}$. We say that $\omega \approx_m \omega'$ if and only if one of the following two compound conditions is satisfied:

- (i) We have $\tau(\omega) = \tau(\omega') = m$ and (ii) $R_m(\omega) = R_m(\omega')$.
- (i) We have $\tau(\omega) = \tau(\omega') < m$, (ii) $R_m(\omega) = R_m(\omega')$ and (iii) the two probability distributions $(\mathbf{P}(G|\mathcal{F}_m)(\omega))_{G \in \tilde{\mathcal{G}}_{m+1}}$ and $(\mathbf{P}(G|\mathcal{F}_m)(\omega'))_{G \in \tilde{\mathcal{G}}_{m+1}}$ have the same support, and their restrictions to their support lie in the same element of $\mathcal{B}_m^c$, where $c = |\{G \in \mathcal{G}_{m+1} : \mathbf{P}(G|\mathcal{F}_m)(\omega) > 0\}|$ is the number of atoms of $\tilde{\mathcal{G}}_{m+1}$ in the common support.

We let $\tilde{\mathcal{G}}_m$ be the partition of $\Omega$ derived by the $\approx_m$ relation.

Define for every $m$ and every atom $G$ of $\tilde{\mathcal{G}}_m$ a probability distribution $\tilde{q}_{G,m}$ over the atoms of $\tilde{\mathcal{G}}_{m+1}$ that are contained in $G$ simply by $\tilde{q}_{G,m} = q_B^*$, where $B$ is the atom of $\mathcal{B}_m^c$ that contains all probability distributions $(\mathbf{P}(G|\mathcal{F}_m)(\omega))_{G \in \tilde{\mathcal{G}}_{m+1}}$ for $\omega \in G$, and $c$ is defined as above.

The sequence $(\tilde{\mathcal{G}}_m)$ satisfies the following properties:

- We have $\tilde{\mathcal{G}}_m \subseteq \mathcal{F}_m$.
- We have that $R_m$ is $\tilde{\mathcal{G}}_m$-measurable.
- For every atom $G$ of $\tilde{\mathcal{G}}_m$, $\sum_{G' \subseteq G, G' \in \tilde{\mathcal{G}}_{m+1}} |\mathbf{P}(G'|\mathcal{F}_m)(\omega) - \tilde{q}_{G,m}(G')| < \delta_m$.

The pair $((\tilde{\mathcal{G}}_m), (\tilde{q}_{G,m}))$ is not a $\delta$ approximation of $\Gamma$ between $n$ and $\tau$, since $(\tilde{\mathcal{G}}_m)$ is not increasing (see condition 4 in Definition 6.2). We remedy this flaw by taking upward intersections. For every $m \geq n$ let $\mathcal{G}_m = \bigvee_{n \leq k \leq m} \tilde{\mathcal{G}}_k$ be the finite $\sigma$-algebra spanned by $\tilde{\mathcal{G}}_n, \ldots, \tilde{\mathcal{G}}_m$. Let $G$ be an atom of $\mathcal{G}_m$. Then $G \subseteq \tilde{G}$ for some atom $\tilde{G}$ of $\tilde{\mathcal{G}}_m$. We define a probability distribution $q_{G,m}$ over the atoms of $\mathcal{G}_{m+1}$ by

$$q_{G,m}(F) = \begin{cases} \tilde{q}_{\tilde{G},m}(G'), & \text{if } F = G \cap G' \text{ for some atom } G' \text{ of } \tilde{\mathcal{G}}_{m+1}, \\ 0, & \text{otherwise.} \end{cases}$$

It can be verified that $((\mathcal{G}_m), (q_{G,m}))$ is a $\delta$ approximation of $\Gamma$ between $n$ and $\tau$, and that each step in the construction is consistent. $\square$

Theorem 6.5, together with the discussion in Section 6.1 that relates $\delta$ approximation to games on trees and Lemmas 6.3 and 6.4, yields the following corollary:

COROLLARY 6.6. *Let $\Gamma$ be a two-player nonzero–sum stopping game. There is a consistent NT function $T = (T_{n,\tau})$ that assigns for every $n \geq 0$, every bounded stopping time $\tau$ and every $\omega \in \Omega$, a game on a tree $T_{n,\tau}(\omega)$ such that the following hold for every $n$ and every $\tau$:*

1. *For every pair of strategies $(x, y)$ in $T_{n,\tau}$ and each $i = 1, 2$, $|\rho^i(x, y; \mathcal{F}, n, \tau) - \rho^i(x, y; T_{n,\tau})| < \Delta_n$.*
2. *For every pair of strategies $(x, y)$ in $T_{n,\tau}$, $|\pi(x, y; \mathcal{F}, n, \tau) - \pi(x, y; T_{n,\tau})| < \Delta_n$.*
3. *If $(x, y)$ is an $\varepsilon$-equilibrium in $T_{n,\tau}(\omega)$ with expected payoff in $[a_1, a_1 + \varepsilon] \times [a_2, a_2 + \varepsilon]$, then:*
   (i) *For every strategy $y'$ of player 2, $\rho^2(x, y'; \mathcal{F}, n, \tau) \leq (a_2 + 2\varepsilon)\pi(x, y'; \mathcal{F}, n, \tau) + \Delta_n$.*
   (ii) *For every strategy $x'$ of player 1, $\rho^1(x', y; \mathcal{F}, n, \tau) \leq (a_1 + 2\varepsilon)\pi(x', y; \mathcal{F}, n, \tau) + \Delta_n$.*

[In the preceding list of inequalities, we identify strategies in $T_{n,\tau}$ with the corresponding $(\mathcal{F}, n, \tau)$-strategies in $\Gamma$.]

**7. Constructing an $\varepsilon$-equilibrium.** In the present section we use all the tools we have developed so far to construct an $\varepsilon$-equilibrium.

7.1. *A sufficient condition for existence of an equilibrium.* Here we provide a sufficient condition for the existence of an equilibrium for games that satisfy the following conditions:

A1. There exists $K \in \mathbf{N}$ such that for every $n \geq 0$, $R_n \in \{0, \pm\frac{1}{K}, \pm\frac{2}{K}, \ldots, \pm\frac{K}{K}\}^6$.
A2. For every $n \geq 0$, $R^1 := \limsup_{n \to \infty} R^1_{\{1\},n}$ is constant and $R^1_{\{1\},n} \leq R^1$.
A3. For every $n \geq 0$, $R^2 := \limsup_{n \to \infty} R^2_{\{2\},n}$ is constant and $R^2_{\{2\},n} \leq R^2$.
A4. Whenever $R^1_{\{1\},n} = R^1$, $R^2_{\{1\},n} < R^2$.
A5. Whenever $R^2_{\{2\},n} = R^2$, $R^1_{\{2\},n} < R^1$.

LEMMA 7.1. *Let $\varepsilon > 0$ be given. Let $\Gamma$ be a stopping game that satisfies conditions* A1–A5. *Suppose that there is a payoff vector $(a_1, a_2)$, an increasing sequence of stopping times $0 = \tau_0 < \tau_1 < \cdots$ and, for every $k \geq 0$, a pair of $(\mathcal{F}, \tau_1, \tau_2)$-strategies $(x_k, y_k)$ such that the following hold for every $k \geq 0$:*

D1. *For $i = 1, 2$, $\rho^i(x_k, y_k; \mathcal{F}, \tau_k, \tau_{k+1}) \geq (a_i - \varepsilon)\pi(x_k, y_k; \mathcal{F}, \tau_k, \tau_{k+1}) - \Delta_{\tau_k}$.*
D2. *For every $(\mathcal{F}, \tau_k, \tau_{k+1})$-strategy $x'_k$, $\rho^1(x'_k, y_k; \mathcal{F}, \tau_k, \tau_{k+1}) \leq (a_1 + \varepsilon)\pi(x'_k, y_k; \mathcal{F}, \tau_k, \tau_{k+1}) + \Delta_{\tau_k}$.*
D3. *For every $(\mathcal{F}, \tau_k, \tau_{k+1})$-strategy $y'_k$, $\rho^2(x_k, y'_k; \mathcal{F}, \tau_k, \tau_{k+1}) \leq (a_2 + \varepsilon)\pi(x_k, y'_k; \mathcal{F}, \tau_k, \tau_{k+1}) + \Delta_{\tau_k}$.*
D4. *Almost surely, $\sum_{k \geq 0} \pi(x_k, y_k; \mathcal{F}, \tau_k, \tau_{k+1}) = \infty$.*
D5. *If $a_1 < R^1 - \varepsilon$, then $\sum_{k \geq 0} \pi(0, y_k; \mathcal{F}, \tau_k, \tau_{k+1}) = \infty$ a.s.*
D6. *If $a_2 < R^2 - \varepsilon$, then $\sum_{k \geq 0} \pi(x_k, 0; \mathcal{F}, \tau_k, \tau_{k+1}) = \infty$ a.s.*

*Then the game admits an $8\varepsilon$-equilibrium.*

In words, conditions D1–D3 roughly say that each $(x_k, y_k)$ is an $\varepsilon$-equilibrium with payoffs close to $(a_1, a_2)$. Condition D4 says that the concatenation of all the $(x_k, y_k)$'s stops with probability 1, while conditions D5 and D6 say that if one player receives a low payoff, then his or her opponent stops with probability 1 under the concatenation. The proof of Lemma 7.1 is standard.

PROOF OF LEMMA 7.1. Denote by $x$ (resp. $y$) the strategy in $\Gamma$ that is the concatenation of the strategies $(x_k)_{k \geq 0}$ [resp. $(y_k)_{k \geq 0}$].

By D4, $\mathbf{P}_{x,y}(\theta < \infty) = 1$. By D5, if $a_1 < R^1 - \varepsilon$, then $\mathbf{P}_{0,y}(\theta < \infty) = 1$. By D6, if $a_2 < R^2 - \varepsilon$, then $\mathbf{P}_{x,0}(\theta < \infty) = 1$. Therefore, there is $L \in \mathbf{N}$ such that

(11) $$\mathbf{P}_{x,y}(\theta < \tau_L) \geq 1 - \varepsilon,$$

(12) $$\begin{aligned} \text{if } a_1 < R^1 - \varepsilon, & \quad \text{then } \mathbf{P}_{0,y}(\theta < \tau_L) \geq 1 - \varepsilon, \\ \text{if } a_2 < R^2 - \varepsilon, & \quad \text{then } \mathbf{P}_{x,0}(\theta < \tau_L) \geq 1 - \varepsilon. \end{aligned}$$

We now define a pair of strategies $(x^*, y^*)$, which is a slight augmentation of $(x, y)$: If $a_1 < R^1 - \varepsilon$, let $y^*$ coincide with $y$. If $a_1 \geq R^1 - \varepsilon$, let $y^*$ be the strategy that coincides with $y$ up to stage $\tau_L$, and from that stage onward, stops with probability $\varepsilon$ whenever $R^2_{\{2\},n} = R^2$ (and with probability 0 otherwise). The strategy $x^*$ is defined analogously.

We argue that $(x^*, y^*)$ is an $8\varepsilon$-equilibrium. We only prove that player 1 cannot profit more than $8\varepsilon$ by deviating.

By summing D1 for $i = 1$ over $k \geq 0$, and since $\mathbf{P}_{x,y}(\theta < \infty) = 1$ and $\sum_{k=1}^{\infty} \Delta_{\tau_k} < \varepsilon$, we obtain

$$\gamma^1(x, y) \geq a_1 - 2\varepsilon. \tag{13}$$

Since $(x^*, y^*)$ coincides with $(x, y)$ up to stage $\tau_L$ and by (11),

$$\gamma^1(x^*, y^*) \geq \gamma^1(x, y) - 2\varepsilon \geq a_1 - 4\varepsilon.$$

Let $x'$ be any strategy of player 1 and, for every $k \geq 0$, let $x'_k$ be the $(\mathcal{F}, \tau_k, \tau_{k+1})$-strategy induced by $x'$.

By summing D2 over $k \geq 0$ and since $\sum_{k=1}^{\infty} \Delta_{\tau_k} < \varepsilon$, we obtain

$$\gamma^1(x', y) \leq (a_1 + \varepsilon)\mathbf{P}_{x',y}(\theta < +\infty) + \varepsilon. \tag{14}$$

If $a_1 < R^1 - \varepsilon$, then by D5, $\mathbf{P}_{x',y}(\theta < +\infty) \geq \mathbf{P}_{0,y}(\theta < +\infty) = 1$. Therefore,

$$\gamma^1(x', y) \leq a_1 + 2\varepsilon.$$

Since $y^*$ coincides with $y$ up to stage $\tau_L$ and by (12),

$$\gamma^1(x', y^*) \leq \gamma^1(x', y) + 2\varepsilon \leq a_1 + 4\varepsilon \leq \gamma^1(x^*, y^*) + 8\varepsilon.$$

We now consider the case $a_1 \geq R^1 - \varepsilon$. By summing D2 over $k = 0, \ldots, L-1$ and since $\sum_{k=1}^{\infty} \Delta_{\tau_k} < \varepsilon$, we obtain

$$\mathbf{E}_{x',y}[R^1_{Q,\theta} \mathbf{1}_{\{\theta < \tau_L\}}] \leq (a_1 + \varepsilon)\mathbf{P}_{x',y}(\theta < \tau_L) + \varepsilon. \tag{15}$$

By the definition of $y^*$ after stage $\tau_L$, it follows that

$$\gamma^1(x', y^*) \leq \mathbf{E}_{x',y}[R^1_{Q,\theta} \mathbf{1}_{\{\theta < \tau_L\}}] + (R^1 + \varepsilon)\mathbf{P}_{x',y}(\theta \geq \tau_L).$$

By (15) and since $R^1 \leq a_1 + \varepsilon$,

$$\gamma^1(x', y^*) \leq a_1 + 2\varepsilon \leq \gamma^1(x^*, y^*) + 6\varepsilon,$$

as desired. □

As the following example shows, adding a threat of punishment might be necessary.

EXAMPLE 7.2. Consider a game with deterministic payoffs: $R_{\{1\},n} = (-1,2)$, $R_{\{2\},n} = (-2,1)$ and $R_{\{1,2\},n} = (0,-3)$. We first argue that all $\varepsilon$-equilibrium payoffs are close to $(-1,2)$.

Given a strategy $x$ of player 1, player 2 can always wait until the probability of stopping under $x$ is exhausted and then stop. Therefore, in any $\varepsilon$-equilibrium, the probability of stopping is at least $1 - \varepsilon$, and the corresponding payoff is close to the convex hull of $(-1,2)$ and $(-2,1)$. Since player 1 can always guarantee $-1$ by stopping at the first stage, the claim follows.

However, in every $\varepsilon$-equilibrium $(x,y)$, we must have $\mathbf{P}_{0,y}(\theta < \infty) \geq 1/2$, otherwise player 1 receives more than $-1$ by never stopping.

Thus, an $\varepsilon$-equilibrium has the following structure, for some integer $N$. Player 1 stops with probability at least $1 - \varepsilon$ before stage $N$ and with probability at most $\varepsilon$ after that stage; player 2 stops with probability at most $\varepsilon$ before stage $N$ and with probability at least $1/2$ after that stage. The strategy of player 2 serves as a threat of punishment: if player 1 does not stop before stage $N$ he or she will be punished in subsequent stages.

7.2. *Proof of Theorem 2.2.* Define $R^1 = \limsup_{n \to \infty} R^1_{\{1\},n}$ and $R^2 = \limsup_{n \to \infty} R^2_{\{2\},n}$. These are the maximal payoffs each player can guarantee by stopping when the opponent always continues. Recall that we assumed that the range $\mathcal{R}$ of the payoff process is finite, so that $p$-a.s. there exist infinitely many $n$'s such that $R^1_{\{1\},n} = R^1$ and infinitely many $n$'s such that $R^2_{\{2\},n} = R^2$.

We fix throughout $\varepsilon > 0$ sufficiently small.

By Corollary 6.6 there is a consistent NT function $(T_{n,\tau})$ that assigns for every $n \geq 0$, every bounded stopping time $\tau$, and every $\omega \in \Omega$, a game on a tree that satisfies several desirable properties.

Define several subsets of $\Omega$ as

$$A_- = \{R^1 \leq 0 \text{ and } R^2 \leq 0\}$$

and, for every $r \in \mathcal{R}$,

$$A_{r,1} = \{(R^1, R^2) = r\} \cap \left\{ \limsup_{\substack{n \to \infty \\ R^1_{\{1\},n} = R^1}} R^2_{\{1\},n} \geq R^2 \right\},$$

$$A_{r,2} = \{(R^1, R^2) = r\} \cap \left\{ \limsup_{\substack{n \to \infty \\ R^2_{\{2\},n} = R^2}} R^1_{\{2\},n} \geq R^1 \right\},$$

$$A_{r,3} = \{(R^1, R^2) = r\} \cap (A_- \cup A_1 \cup A_2)^c.$$

Recall that for every two payoff processes $R$ and $\widetilde{R}$ such that $\mathbf{E}[\sup_{n\in\mathbf{N}}\|R_n - \widetilde{R}_n\|_\infty] < \varepsilon$, every $\varepsilon$-equilibrium in $(\Omega, \mathcal{A}, p, \mathcal{F}, \widetilde{R})$ is a $3\varepsilon$-equilibrium in $(\Omega, \mathcal{A}, p, \mathcal{F}, R)$.

We are now going to further partition the sets $(A_{r,3})_{r\in\mathcal{R}}$. In Section 5.3 we attached to each game on a graph $T$ a color $c(T)$ (from a finite set of possible colors) and $J$ numbers in the unit interval $(\lambda_j(T))_{j=1}^J$.

By Theorem 4.3 there is an increasing sequence of stopping times $0 \leq \tau_0 < \tau_1 < \cdots$ such that

$$\mathbf{P}(c(T_{\tau_0,\tau_1}) = c(T_{\tau_k,\tau_{k+1}}) \ \forall k \in \mathbf{N}) > 1 - \varepsilon.$$

Define the subsets of $A_{r,3}$, for every $j = 1, \ldots, J$,

$$A_{r,3;j} = A_{r,3} \cap \left\{\sum_{k\geq 0} \lambda_j(T_{\tau_k,\tau_{k+1}}) = \infty\right\}$$

and, for each good rectangle $g \in G$,

$$B_{r,3;g} = A_{r,3} \cap \{c(T_{\tau_k,\tau_{k+1}}) = g \ \forall k\}.$$

Let $\mathcal{P}^*$ be the collection of sets $A_-$, $(A_{r,1})$, $(A_{r,2})$, $(A_{r,3;j})$, $(B_{r,3;g})$. By properly modifying the payoff process on a set with measure smaller than $\varepsilon$, and by possibly dropping the first few stopping times in the sequence $(\tau_k)$ and renaming the remaining stopping times, we can assume w.l.o.g. that (i) all the sets in $\mathcal{P}^*$ are $\mathcal{F}_{\tau_0}$-measurable, (ii) $R^1_{\{1\},k} \leq R^1$ and $R^1_{\{2\},k} \leq R^2$ for every $k \geq \tau_0$, and (iii) $\sum_{j=1}^J \sum_{k\geq 0} \lambda_j(T_{\tau_k,\tau_{k+1}}) \leq \varepsilon$ on $\bigcup_{g\in G} B_{r,3;g}$.

As the sets in $\mathcal{P}^*$ are not necessarily disjoint, we let $\mathcal{P}$ be a finite partition of $\Omega$ into $\mathcal{F}_{\tau_0}$-measurable sets such that for every atom $B$ of $\mathcal{P}$ and every set $B^*$ in $\mathcal{P}^*$, either $B \subseteq B^*$ or $B \cap B^* = \varnothing$.

For every $\mathcal{F}_{\tau_0}$-measurable set $B$ define the stopping game that is restricted to $B$ and starts at stage $\tau_0$: $\Gamma_{B,\tau_0} = (B, \mathcal{A}_B, p_{|B}, (\mathcal{F}_{\tau_0+k})_{k\geq 0}, (R_{\tau_0+k})_{k\geq 0})$, where $\mathcal{A}_B$ is the $\sigma$-algebra over $B$ induced by $\mathcal{A}$ and $p_{|B}$ is the probability distribution $p$ conditioned on $B$.

The following lemma is standard.

LEMMA 7.3. *If for every atom $B$ of $\mathcal{P}$ with $\mathbf{P}(B) > 0$, the game $\Gamma_{B,\tau_0}$ admits an $\varepsilon$-equilibrium, then the game $\Gamma$ admits a $3\varepsilon$-equilibrium.*

PROOF. It is well known that any finite-stage game admits a 0-equilibrium [one can, e.g., adapt the proof for the zero–sum case given in Rosenberg, Solan and Vieille (2001), Proposition 3.1]. Since $\tau_0$ is bounded, and since we changed payoffs on a set of measure at most $\varepsilon$, the following strategy pair is a $3\varepsilon$-equilibrium:

- If the game has not terminated by stage $\tau_0$, play from that stage on an $\varepsilon$-equilibrium in $\Gamma_{B,\tau_0}$, where $B$ is the atom of $\mathcal{P}$ that contains $\omega$. Denote by $\gamma$ the $\mathcal{F}_{\tau_0}$-measurable payoff that corresponds to this strategy.

- Until state $\tau_0$, play a 0-equilibrium in the game that terminates at stage $\tau_0$, if no player stops before that stage, with terminal payoff $\gamma$.  □

Thus, it remains to show that for every atom $B$ of $\mathcal{P}$, the game $\Gamma_{B,\tau_0}$ admits a $C\varepsilon$-equilibrium for some $C > 0$. We verify this in the next sections. For convenience, we assume w.l.o.g. that $\tau_0 = 0$.

7.2.1. $B \subseteq A_-$. The game $\Gamma_{B,n}$ admits a 0-equilibrium: The strategy pair where no player ever stops is a 0-equilibrium. Indeed, under this strategy the expected payoff for both players is 0, while any player $i$ who stops at some stage $k$ receives $R^i_{\{i\},k} \leq R^i \leq 0$.

7.2.2. (i) $B \subseteq A_{r,1}$ and $r^1 \geq 0$ or (ii) $B \subseteq A_{r,2}$ and $r^2 \geq 0$. We deal only with case (i). The following strategy pair is a $2\varepsilon$-equilibrium: Player 2 never stops, while player 1 stops with probability $\varepsilon$ at every stage $k$ such that $R^1_{\{1\},k} = R^1 = r^1$ and $R^2_{\{1\},k} \geq R^2$. The expected payoff to player 1 is $r^1$, while the expected payoff to player 2 is at least $r^2$. Since $r^1 \geq 0$, player 1 cannot profit by deviating, while if player 2 stops at some stage $k$, his or her expected payoff is bounded by $(1-\varepsilon)R^2_{\{2\},k} + \varepsilon \leq (1-\varepsilon)r^2 + \varepsilon \leq r^2 + 2\varepsilon$.

7.2.3. (i) $B \subseteq A_{r,1}$ and $r^1 < 0$ or (ii) $B \subseteq A_{r,2}$ and $r^2 < 0$. We deal only with case (i). If $r^2 \leq 0$, then $B \subseteq A_-$, so that by Section 7.2.1 $\Gamma_{B,n}$ admits a 0-equilibrium. Assume then that $r^2 > 0$. If $\limsup_{n \to \infty, R^2_{\{2\},n} = R^2} R^1_{\{2\},n} \geq R^1$, then $B \subseteq A_{r,\{2\}}$, so that by Section 7.2.2 the game $\Gamma_{B,n}$ admits a $2\varepsilon$-equilibrium. Assume then that $\limsup_{n \to \infty, R^2_{\{2\},n} = R^2} R^1_{\{2\},n} < R^1$. Under these assumptions, the following strategy pair is a $3\varepsilon$-equilibrium: Player 1 stops at every stage $k$ such that $R^1_{\{1\},k} = r^1$ and $R^2_{\{1\},k} \geq r^2$ with probability $\varepsilon$, while player 2 stops at every stage $k > N$ such that $R^1_{\{2\},k} < r^1$ with probability $\varepsilon$, where $N$ is sufficiently large so that under this strategy, player 1 stops with probability at least $1 - \varepsilon$ prior to stage $N$.

In the last two cases we use the coloring procedure presented in Section 5.3 and the sufficient condition that appears in Lemma 7.1.

7.2.4. $B \subseteq A_{r,3;j}$. Recall that on $A_{r,3;j}$,

$$\sum_{k \geq 0} \lambda_j(T_{\tau_k, \tau_{k+1}}) = \infty. \tag{16}$$

Denote by $T^{(j)}(\tau_k, \tau_{k+1})$ the tree that was generated in round $j$ of the coloring procedure in Section 5.3 and denote by $(x^{(j)}(\tau_k, \tau_{k+1}), y^{(j)}(\tau_k, \tau_{k+1}))$

NONZERO–SUM STOPPING GAMES 23the $\varepsilon$-equilibrium with payoff in a bad rectangle $[a_1, a_1 + \varepsilon] \times [a_2, a_2 + \varepsilon]$. By the definition of $\lambda_j$,

$$(17) \quad \pi(x^{(j)}(\tau_k, \tau_{k+1}), y^{(j)}(\tau_k, \tau_{k+1}); T_{\tau_k, \tau_{k+1}}) \geq \varepsilon^2 \times \lambda_j(T_{\tau_k, \tau_{k+1}}).$$

We apply Lemma 7.1 with $2\varepsilon$, $\{\tau_k\}$ and $\{(x^{(j)}(\tau_k, \tau_{k+1}), y^{(j)}(\tau_k, \tau_{k+1}))\}$. Since the rectangle is bad, $a_1 \geq R^1 - \varepsilon$ and $a_2 \geq R^2 - \varepsilon$, so that conditions D5 and D6 trivially hold.

Condition D4 holds by (16), (17) and Corollary 6.6.2. Condition D1 follows from the choice of $(x^{(j)}(\tau_k, \tau_{k+1}), y^{(j)}(\tau_k, \tau_{k+1}))$ and by Corollary 6.6. Condition D2 follows from Corollary 6.6.3. The proof that condition D3 holds is analogous. The desired result follows by Lemma 7.1.

7.2.5. $B \subseteq B_{r,3;g}$. Recall that on $B_{r,3;g}$, $c(T(\tau_k, \tau_{k+1})) = g$ for every $k$ and $\sum_{j=1}^{J} \sum_{k \geq 0} \lambda_j(T_{\tau_k, \tau_{k+1}}) \leq \varepsilon$. Since $g = [a_1, a_1 + \varepsilon] \times [a_2, a_2 + \varepsilon]$ is a good rectangle, $a_1 \leq R^1 - 2\varepsilon$ or $a_2 \leq R^2 - 2\varepsilon$. We assume w.l.o.g. that $a_1 \leq R^1 - 2\varepsilon$. Let $(x^{(0)}(\tau_k, \tau_{k+1}), y^{(0)}(\tau_k, \tau_{k+1}))$ be an $\frac{\varepsilon}{2}$-equilibrium in $T^{(J)}(\tau_k, \tau_{k+1})$ with payoff in $g$. Since $a_1 \leq R^1 - 2\varepsilon$, Lemma 5.3 implies that

$$(18) \quad \begin{aligned} \pi(0, y^{(0)}(\tau_k, \tau_{k+1}); T_{\tau_k, \tau_{k+1}}) &\geq \frac{\varepsilon}{6} \mu_1(T^{(J)}(\tau_k, \tau_{k+1})) \\ &\geq \frac{\varepsilon}{6} \mu_1(T(\tau_k, \tau_{k+1})) - \sum_{j=1}^{J} \lambda_j(T_{\tau_k, \tau_{k+1}}). \end{aligned}$$

We apply Lemma 7.1, with $2\varepsilon$, $\{\tau_k\}$ and $\{(x^{(0)}(\tau_k, \tau_{k+1}), y^{(0)}(\tau_k, \tau_{k+1}))\}$. Condition D1 holds by the choice of $(x^{(0)}(\tau_k, \tau_{k+1}), y^{(0)}(\tau_k, \tau_{k+1}))$ and by Corollary 6.6. Conditions D2 and D3 hold as in the previous case.

We now prove that condition D5 holds. By (18) and Corollary 6.6.2,

$$\pi(0, y^{(0)}(\tau_k, \tau_{k+1}); \mathcal{F}, \tau_k, \tau_{k+1}) \geq \frac{\varepsilon}{6} \mu_1(T(\tau_k, \tau_{k+1})) - \sum_{j=1}^{J} \lambda_j(T_{\tau_k, \tau_{k+1}}) - \Delta_{\tau_k}.$$

Since $R^1_{\{1\},n} = R^1$ i.o. and since $\sum_{j=1}^{J} \sum_{k \geq 0} \lambda_j(T_{\tau_k, \tau_{k+1}}) \leq \varepsilon$ on $B_{r,3;g}$, it follows that condition D5 holds.

Condition D6 follows by an analog argument. Condition D4 follows from conditions D5 and D6, and since $a_1 \leq R^1 - 2\varepsilon$ or $a_2 \leq R^2 - 2\varepsilon$.

**8. Proof of Proposition 5.5.** This section is devoted to the proof of Proposition 5.5. We start by defining some new concepts for games on finite trees and by studying some of their properties. The proof itself appears in Section 8.5.

Throughout the section we fix a stopping game on a tree that satisfies conditions B1–B4.

8.1. *Union of strategies.* Given $n$ stationary strategies $x_1, x_2, \ldots, x_n$, we define their *union* $x$ by $x(s) = 1 - \prod_{1 \leq k \leq n}(1 - x_k(s))$. The probability that the union strategy continues at each node is the probability that *all* of its components continue. We denote $x = x_1 \dotplus x_2 \dotplus \cdots \dotplus x_n$. Given $n$ pairs of stationary strategies $\alpha_k = (x_k, y_k)$, $1 \leq k \leq n$, we denote by $\alpha_1 \dotplus \cdots \dotplus \alpha_n$ the stationary strategy pair $(x, y)$ that is defined by $x = x_1 \dotplus \cdots \dotplus x_n$ and $y = y_1 \dotplus \cdots \dotplus y_n$.

Consider $n$ copies of the game that are played simultaneously, such that the choice of a new node is the same across the copies; that is, all copies that have not terminated at stage $t$ are at the same node. Nevertheless, the lotteries made by the players concerning the decision whether to stop are independent. Let $\alpha_k = (x_k, y_k)$, $1 \leq k \leq n$, be the stationary strategy pair used in copy $k$ and let $\alpha = \alpha_1 \dotplus \cdots \dotplus \alpha_n$.

We consider the first round of the game. Let $t_k$ be the stopping stage in copy $k$, let $s_k$ be the node in which termination occurred, let $Q_k$ be the set of players that stop at stage $t_k$ and let $r_k^i = R_{Q_k,s_k}^i \mathbf{1}_{\{t_k < \infty\}}$ be the payoff in the first round. Set

$$\pi_k = \pi(x_k, y_k) = \mathbf{P}(t_k < \infty) \quad \text{and} \quad \rho_k^i = \rho^i(x_k, y_k) = \mathbf{E}[r_k^i].$$

Let $t$, $r$, $\rho$ and $\pi$ be the analog quantities w.r.t. $\alpha$: Denoting $k^*$ the index of a copy that stops first (so that $t_{k^*} = \min_{k=1,\ldots,n} t_k$), we have $t = t_{k^*}$, $s = s_{k^*}$, $Q = \bigcup_{k:t_k=t_{k^*}} Q_k$ and $r^i = R_{Q,s}^i \mathbf{1}_{t<+\infty}$. Moreover, $\rho^i = \rho^i(x, y) = \mathbf{E}[r^i]$ and $\pi = \pi(x, y) = \mathbf{P}(t < \infty)$.

Let $\gamma_k = \gamma(x_k, y_k)$ be the expected payoff under $\alpha_k = (x_k, y_k)$ and let $\gamma = \gamma(x, y)$ be the corresponding quantity under $\alpha$.

The following lemma follows from the independence of the plays given the branch. Recall that $F_s$ is the event that the randomly chosen branch passes through $s$.

LEMMA 8.1. *Let $s \in S_0$ be a node of depth $j$. Then, for every $1 \leq k, l \leq n$, $l \neq k$, the event $\{t_k \leq j\}$ and the random variable $t_k \mathbf{1}_{\{t_k \leq j\}}$ are independent of $t_l$ given $F_s$.*

LEMMA 8.2. *Let $N = \sum_{k=1}^n \mathbf{1}_{\{t_k < \infty\}}$ be the number of copies that terminate in the first round. Then:*

1. *We have $\sum_{k=1}^n \pi_k - \mathbf{E}[N\mathbf{1}_{\{N \geq 2\}}] \leq \pi \leq \sum_{k=1}^n \pi_k$.*
2. *We have $\sum_{k=1}^n \rho_k^i - \mathbf{E}[(N+1)\mathbf{1}_{\{N \geq 2\}}] \leq \rho^i \leq \sum_{k=1}^n \rho_k^i + \mathbf{E}[(N+1)\mathbf{1}_{\{N \geq 2\}}]$ for each player $i \in \{1, 2\}$.*

PROOF. Observe that

$$N - N\mathbf{1}_{\{N \geq 2\}} = \mathbf{1}_{\{N=1\}} \leq \mathbf{1}_{\{N \geq 1\}} \leq N = \sum_{k=1}^n \mathbf{1}_{\{t_k < \infty\}}.$$

The first result follows by taking expectations.

For the second result, note that

$$\text{(19)} \quad \sum_{k=1}^{n} r_k^i - (N+1)\mathbf{1}_{\{N \geq 2\}} \leq r^i \leq \sum_{k=1}^{n} r_k^i + (N+1)\mathbf{1}_{\{N \geq 2\}}.$$

Indeed, on $\{N \leq 1\}$, (19) holds with equality and on $\{N \geq 2\}$, the left-hand side is at most $-1$, whereas the right-hand side is at least $+1$. The result follows by taking expectations. $\square$

### 8.2. Heavy and light nodes.

DEFINITION 8.3. Let $\sigma = (x, y)$ be a pair of stationary strategies and let $\delta > 0$. A node $s \in S_0$ is $\delta$-*heavy* with respect to $\sigma$ if $\mathbf{P}_\sigma(t < \infty | F_s) \geq \delta$; that is, the probability of termination in the first round given that the chosen branch passes through $s$ is at least $\delta$. The node $s$ is $\delta$-*light* w.r.t. $\sigma$ if $\mathbf{P}_\sigma(t < \infty | F_s) < \delta$.

For a fixed $\delta$, we denote by $H_\delta(\sigma)$ the set of $\delta$-heavy nodes w.r.t. $\sigma$. Two simple implications of this definition follow:

FACT 1. We have $H_\delta(\alpha_1) \subseteq H_\delta(\alpha_1 \dotplus \alpha_2)$.

FACT 2. We have $H_{\delta_1}(\sigma) \subseteq H_{\delta_2}(\sigma)$ whenever $\delta_1 \geq \delta_2$.

The following lemma asserts that for every $\varepsilon$-equilibrium $(x, y)$ with high payoffs, $H_\varepsilon(x, y)$ is nonempty, regardless of the size of the tree.

LEMMA 8.4. *Let* $\varepsilon \in (0, 1/(36K)^2)$ *and let* $(x, y)$ *be a stationary* $\varepsilon$-*equilibrium such that* $R^i - \varepsilon \leq \gamma^i(x, y)$, $i = 1, 2$. *Then* $H_\varepsilon(x, y) \neq \phi$. *In particular, by Fact 2,* $H_{\varepsilon^2}(x, y) \neq \phi$.

COMMENT. The proof hinges on the assumption that $R^2_{\{1\},s} < R^2$ whenever $R^1_{\{1\},s} = R^1$. As a counterexample when this condition does not hold, take a game in which (a) $R^i_{Q,s} = 1$ for every $i$, $Q$ and $s$, and (b) $R^1 = R^2 = 1$. Then any stationary strategy pair which stops with positive probability is a 0-equilibrium.

PROOF OF LEMMA 8.4. We prove that there is a node $s$ such that $x_s y_s > \varepsilon$, so that $s \in H_\varepsilon(x, y)$. The idea is that if $H_\varepsilon(x, y)$ is empty, then the probability that both players stop simultaneously under $(x, y)$ is small, but if $(x, y)$ is an $\varepsilon$-equilibrium, this probability must be high.

Denote $\pi_1 = \mathbf{P}_{x,y}(t < \infty, Q = \{1\})$, $\pi_2 = \mathbf{P}_{x,y}(t < \infty, Q = \{2\})$ and $\pi_{1,2} = \mathbf{P}_{x,y}(t < \infty, Q = \{1, 2\})$. These are the probabilities that in a single round, player 1 stops alone, player 2 stops alone and both players stop simultaneously, respectively.

Since for some player $i$, $R^i$ is strictly positive, and hence strictly larger than $\varepsilon$, we have $\gamma^i(x,y) \geq R^i - \varepsilon > 0$. Therefore, $\pi_1 + \pi_2 + \pi_{1,2} > 0$. Assume w.l.o.g. that $\pi_1 \geq \pi_2$.

Suppose to the contrary that $x_s y_s \leq \varepsilon$ for every node $s \in S_0$. The probability that both players stop simultaneously at $s$, conditional on the game terminating at $s$, is
$$f(x_s, y_s) = \frac{x_s y_s}{x_s(1 - y_s) + (1 - x_s)y_s + x_s y_s}.$$

The maximum of the function $f(x_s, y_s)$ over the set $\{x_s y_s \leq \varepsilon, 0 \leq x_s, y_s \leq 1\}$ is attained at the point $x_s = y_s = \sqrt{\varepsilon}$. Therefore, $f(x_s, y_s) \leq f(\sqrt{\varepsilon}, \sqrt{\varepsilon}) \leq \sqrt{\varepsilon}$. This implies that $\pi_{1,2} \leq \sqrt{\varepsilon}$.

Since $(x, y)$ is an $\varepsilon$-equilibrium and since payoffs are bounded by 1,
$$R^2 - \varepsilon \leq \gamma^2(x, y)$$
$$\leq \frac{\pi_1 \max_{s \in S_0} R^2_{\{1\}, s} + \pi_2 \max_{s \in S_0} R^2_{\{2\}, s} + \pi_{1,2} \max_{s \in S_0} R^2_{\{1,2\}, s}}{\pi_1 + \pi_2 + \pi_{1,2}}$$
$$\leq \frac{\pi_1(R^2 - 1/K) + \pi_2 R^2 + \pi_{1,2}}{\pi_1 + \pi_2 + \pi_{1,2}}.$$

Since $\pi_1 \geq \pi_2$ and $\varepsilon < 1/4K$, this implies that $\pi_{1,2} \geq 1/6K$. In particular, $\sqrt{\varepsilon} \geq \pi_{1,2} \geq 1/6K$, which is a contradiction when $\varepsilon < 1/36K^2$. □

### 8.3. Orthogonal strategies.

DEFINITION 8.5. Let $\delta > 0$. A sequence $(\alpha_1, \alpha_2, \ldots, \alpha_n)$ of stationary strategy pairs is $\delta$-*orthogonal* if $\alpha_{k+1}(s) = (0, 0)$ for every $1 \leq k \leq n - 1$ and every node $s \in H_\delta(\alpha_1 \dotplus \cdots \dotplus \alpha_k)$; that is, $\alpha_{k+1}$ continues on $\delta$-heavy nodes of $\alpha_1 \dotplus \cdots \dotplus \alpha_k$.

LEMMA 8.6. *Let $\delta > 0$, let $(\alpha_1, \ldots, \alpha_n)$ be a $\delta$-orthogonal sequence of stationary strategy pairs, let $k \in \{1, \ldots, n\}$ and let $s \in S$ be a node of depth $j$. Then*

(20) $\quad \mathbf{P}\left(\{j \leq t_k < \infty\} \cap \left(\bigcup_{l < k} \{t_l < \infty\}\right) \Big| F_s\right) \leq \delta \times \mathbf{P}(j \leq t_k < \infty | F_s).$

PROOF. Fix $k \in \{1, \ldots, n\}$. We prove the lemma by induction on the nodes of $T$, starting from the leaves and climbing up to the root.

Let $s \in S_1$ be a leaf of $T$. Since $s$ is a leaf, $\mathbf{P}(j \leq t_k < \infty) = 0$ and (20) is trivially satisfied.

Assume now that $s \in S_0$. Then

$$\mathbf{P}\left(\{j \leq t_k < \infty\} \cap \left(\bigcup_{l<k}\{t_l < \infty\}\right)\Big|F_s\right)$$

(21)
$$= \mathbf{P}\left(\{t_k = j\} \cap \left(\bigcup_{l<k}\{t_l < \infty\}\right)\Big|F_s\right)$$

$$+ \sum_{s' \in C_s} p_s[s'] \times \mathbf{P}\left(\{j+1 \leq t_k < \infty\} \cap \left(\bigcup_{l<k}\{t_l < \infty\}\right)\Big|F_{s'}\right).$$

By the induction hypothesis, for every child $s' \in C_s$,

$$\mathbf{P}\left(\{j+1 \leq t_k < \infty\} \cap \left(\bigcup_{l<k}\{t_l < \infty\}\right)\Big|F_{s'}\right) \leq \delta \times \mathbf{P}(j+1 \leq t_k < \infty|F_{s'}).$$

(22)

By Lemma 8.1, $\{t_k = j\}$ and $\bigcup_{l<k}\{t_l < \infty\}$ are independent given $F_s$. Therefore

$$\mathbf{P}\left(\{t_k = j\} \cap \left(\bigcup_{l<k}\{t_l < \infty\}\right)\Big|F_s\right) = \mathbf{P}(t_k = j|F_s) \times \mathbf{P}\left(\bigcup_{l<k}\{t_l < \infty\}\Big|F_s\right).$$

If $s$ is $\delta$-light w.r.t. $\alpha_1 \dotplus \cdots \dotplus \alpha_{k-1}$, then $\mathbf{P}(\bigcup_{l<k}\{t_l < \infty\}|F_s) < \delta$. If $s$ is $\delta$-heavy, then, according to the definition of orthogonality, $\mathbf{P}(t_k = j|F_s) = 0$. Therefore,

(23) $$\mathbf{P}\left(\{t_k = j\} \cap \left(\bigcup_{l<k}\{t_l < \infty\}\right)\Big|F_s\right) \leq \delta \times \mathbf{P}(t_k = j|F_s).$$

Equations (21)–(23) yield

$$\mathbf{P}\left(\{j \leq t_k < \infty\} \cap \left(\bigcup_{l<k}\{t_l < \infty\}\right)\right)$$

$$\leq \delta \times \mathbf{P}(t_k = j|F_s) + \delta \times \sum_{s' \in C_s} p_s[s'] \times \mathbf{P}(j+1 \leq t_k < \infty|F_{s'})$$

$$= \delta \times \mathbf{P}(j \leq t_k < \infty|F_s),$$

as desired. □

Applying Lemma 8.6 to the root we get:

COROLLARY 8.7. *Let $\delta > 0$ and let $(\alpha_1, \ldots, \alpha_n)$ be a $\delta$-orthogonal sequence of stationary strategy pairs. For every $k \in \{1, \ldots, n\}$,*

$$\mathbf{P}\left(\{t_k < \infty\} \cap \left(\bigcup_{l<k}\{t_l < \infty\}\right)\right) \leq \delta \times \mathbf{P}(\{t_k < \infty\}) = \delta \pi_k.$$

LEMMA 8.8. *Let $\delta > 0$, let $(\alpha_1, \ldots, \alpha_n)$ be a $\delta$-orthogonal sequence of stationary strategy pairs and let $N = \sum_{k=1}^n \mathbf{1}_{\{t_k < \infty\}}$. Then $\mathbf{E}[(N+1)\mathbf{1}_{\{N \geq 2\}}] \leq 3\delta(\pi_1 + \pi_2 + \cdots + \pi_n)$.*

PROOF. Observe that $N + 1 \leq 3(N-1)$ on $\{N \geq 2\}$ and $(N-1)\mathbf{1}_{\{N \geq 2\}} = \sum_{k=1}^n \mathbf{1}_{\{t_k<\infty\} \cap (\bigcup_{l<k}\{t_l<\infty\})}$. Therefore,

$$\mathbf{E}[(N+1)\mathbf{1}_{\{N\geq 2\}}] \leq 3\mathbf{E}[(N-1)\mathbf{1}_{\{N\geq 2\}}]$$
$$= 3\sum_{k=1}^n \mathbf{P}\left(\{t_k < \infty\} \cap \left(\bigcup_{l<k}\{t_l < \infty\}\right)\right).$$

The result follows by Corollary 8.7. □

From Lemmas 8.2 and 8.8 we get the following:

COROLLARY 8.9. *Let $\delta > 0$ and let $(\alpha_1, \ldots, \alpha_n)$ be a $\delta$-orthogonal sequence of strategy pairs. Denote $\alpha = \alpha_1 \dotplus \cdots \dotplus \alpha_n$. Then for $i = 1, 2$:*

1. *We have $(1 - 3\delta)\sum_{k=1}^n \pi_k \leq \pi \leq \sum_{k=1}^n \pi_k$.*
2. *We have $\sum_{k=1}^n \rho_k^i - 3\delta \sum_{k=1}^n \pi_k \leq \rho^i \leq \sum_{k=1}^n \rho_k^i + 3\delta \sum_{k=1}^n \pi_k$.*

LEMMA 8.10. *Let $\delta > 0$ and let $(\alpha_1, \ldots, \alpha_n)$ be a $\delta$-orthogonal sequence of stationary strategy pairs. Denote $\alpha = \alpha_1 \dotplus \cdots \dotplus \alpha_n$. Then for $i = 1, 2$,*

$$\sum_{k=1}^n \rho_k^i - 6\delta \sum_{k=1}^n \pi_k \leq \gamma^i \times \sum_{k=1}^n \pi_k \leq \sum_{k=1}^n \rho_k^i + 6\delta \sum_{k=1}^n \pi_k.$$

PROOF. By Corollary 8.9 and (6),

$$\sum_{k=1}^n \rho_k^i - 3\delta \sum_{k=1}^n \pi_k \leq \rho^i = \gamma^i \times \pi \leq \begin{cases} \gamma^i \times \sum_{k=1}^n \pi_k, & \text{if } \gamma^i > 0, \\ \gamma^i(1 - 3\delta)\sum_{k=1}^n \pi_k, & \text{if } -1 \leq \gamma^i \leq 0. \end{cases}$$

In both cases, the right-hand side is bounded by $\gamma^i \times \sum_{k=1}^n \pi_k + 3\delta \sum_{k=1}^n \pi_k$, so that

$$\sum_{k=1}^n \rho_k^i - 6\delta \sum_{k=1}^n \pi_k \leq \gamma^i \times \sum_{k=1}^n \pi_k.$$

The proof of the right-hand inequality is similar. □

From Lemma 8.10 and (6) we get:

COROLLARY 8.11. *Let $\delta > 0$ and let $(\alpha_1, \ldots, \alpha_n)$ be a $\delta$-orthogonal sequence of stationary strategy pairs. Denote $\alpha = \alpha_1 \dot{+} \cdots \dot{+} \alpha_n$. Let $-1 \leq u, v \leq 1$.*

1. *If $u \leq \gamma_k^i$ for each $k \in \{1, \ldots, n\}$, then $u - 6\delta \leq \gamma^i$.*
2. *If $\gamma_k^i \leq v$ for each $k \in \{1, \ldots, n\}$, then $\gamma^i \leq v + 6\delta$.*

8.4. *Strong orthogonality.* In the present section we define a stronger notion of orthogonality and study its properties.

DEFINITION 8.12. Let $\delta > 0$. A sequence $(\alpha_1, \alpha_2, \ldots, \alpha_n)$ of stationary strategy pairs is $\delta$-*strongly orthogonal* if, for every $k \in \{1, \ldots, n-1\}$ and every node $s \in H_\delta(\alpha_1 \dot{+} \cdots \dot{+} \alpha_k)$, $\alpha_{k+1}(s') = (0,0)$ for $s' = s$ and for every descendent $s'$ of $s$; that is, $\alpha_{k+1}$ continues from $s$ onward.

The following lemma provides a way to construct $\varepsilon$-orthogonal sequences of strategy pairs from a single $\varepsilon^2$-strongly orthogonal sequence.

LEMMA 8.13. *Let $\varepsilon > 0$ and let $y_1, y_2, \ldots, y_n$ be stationary strategies of player 2 such that the sequence $((0, y_1), \ldots, (0, y_n))$ is $\varepsilon^2$-strongly orthogonal. Let $\bar{x}$ be any pure stationary strategy of player 1 that does not stop twice on the same branch; that is, if $\bar{x}(s) = 1$, then $\bar{x}(s') = 0$ for every descendent $s'$ of $s$. Define strategies $(\bar{x}_k)_{k=1}^n$ of player 1 in the following way: For each $s \in S$ such that $\bar{x}(s) = 1$, let $\bar{x}_k(s) = 1$, where $k \leq n$ is the greatest index for which $s \notin H_\varepsilon((0, y_1) \dot{+} \cdots \dot{+} (0, y_{k-1}))$. Define $\bar{x}_k(s) = 0$ otherwise. Let $\bar{\alpha}_k = (\bar{x}_k, y_k)$. Then the sequence $(\bar{\alpha}_1, \ldots, \bar{\alpha}_n)$ is $\varepsilon$-orthogonal.*

PROOF. By the definition of $(\bar{x}_k)_{1 \leq k \leq n}$ and Fact 1, we get, for every $l \in \{1, \ldots, n-1\}$:

(24)
$$\text{If } \bar{x}_l(s) = 1 \quad \text{then } s \in H_\varepsilon((0, y_1) \dot{+} \cdots \dot{+} (0, y_l));$$
$$\text{If } \bar{x}_{l+1}(s) = 1 \quad \text{then } s \notin H_\varepsilon((0, y_1) \dot{+} \cdots \dot{+} (0, y_l)).$$

Let $l \in \{1, \ldots, n-1\}$ and let $s \in S$ be $\varepsilon$-heavy with respect to $\bar{\sigma}_l = \bar{\alpha}_1 \dot{+} \cdots \dot{+} \bar{\alpha}_l$. We prove that $\bar{x}_{l+1}(s) = y_{l+1}(s) = 0$.

We first prove that $\bar{x}_{l+1}(s) = 0$. Since $s$ is $\varepsilon$-heavy w.r.t. $\bar{\sigma}_l = \bar{\alpha}_1 \dot{+} \cdots \dot{+} \bar{\alpha}_l$, $\mathbf{P}_{\bar{\sigma}_l}(t < \infty | F_s) \geq \varepsilon$. Assume to the contrary that $\bar{x}_{l+1}(s) = 1$. By (24), $s$ is $\varepsilon$-light w.r.t. $(0, y_1) \dot{+} \cdots \dot{+} (0, y_l)$ and, therefore, $\mathbf{P}_{(0,y_1) \dot{+} \cdots \dot{+} (0,y_l)}(t < \infty | F_s) < \varepsilon$. It follows that $\mathbf{P}_{(\bar{x}_1,0) \dot{+} \cdots \dot{+} (\bar{x}_l,0)}(t < \infty | F_s) > 0$, a contradiction to the assumption that $\bar{x}$ does not stop twice on the same branch.

We proceed to prove that $y_{l+1}(s) = 0$. Assume first that there exists an ancestor $s'$ of $s$ such that $\bar{x}_1(s') + \cdots + \bar{x}_l(s') = 1$. By (24) and Fact 1, $s' \in H_\varepsilon((0, y_1) \dotplus \cdots \dotplus (0, y_l))$. Since $((0, y_1), \ldots, (0, y_n))$ is $\varepsilon$-strongly orthogonal, $y_{l+1}(s) = 0$.

We assume now that $\bar{x}_1(s') + \cdots + \bar{x}_l(s') = 0$ for every ancestor $s'$ of $s$. Let $\tilde{D}$ be the (possibly empty) set of $s$'s descendants $d$ that are $\varepsilon$-heavy w.r.t. $(0, y_1) \dotplus \cdots \dotplus (0, y_l)$ and let $D$ be the set that is obtained by removing from $\tilde{D}$ all nodes that have strict ancestors in $\tilde{D}$. By the definition of $D$, $\mathbf{P}_{(0,y_1) \dotplus \cdots \dotplus (0,y_l)}(t < \infty | F_d) \geq \varepsilon$ for every $d \in D$. Let $Y = \bigcup_{d \in D} F_d$. Since this is a mutually disjoint union, it follows that if $Y \neq \phi$, then

$$\mathbf{P}_{(0,y_1) \dotplus \cdots \dotplus (0,y_l)}(t < \infty | Y) \geq \varepsilon \geq \varepsilon \times \mathbf{P}_{\bar{\sigma}_l}(t < \infty | Y).$$

By (24) and the definition of $(\bar{x}_k)_{1 \leq k \leq n}$, it follows that

$$\mathbf{P}_{(0,y_1) \dotplus \cdots \dotplus (0,y_l)}(t < \infty | Y^c \cap F_s) = \mathbf{P}_{\bar{\sigma}_l}(t < \infty | Y^c \cap F_s)$$
$$\geq \varepsilon \times \mathbf{P}_{\bar{\sigma}_l}(t < \infty | Y^c \cap F_s).$$

Combining the last two inequalities and observing that $Y \subseteq F_s$, we get

$$\mathbf{P}_{(0,y_1) \dotplus \cdots \dotplus (0,y_l)}(t < \infty | F_s) \geq \varepsilon \times \mathbf{P}_{\bar{\sigma}_l}(t < \infty | F_s) \geq \varepsilon^2.$$

Thus $s$ is $\varepsilon^2$-heavy with respect to $(0, y_1) \dotplus \cdots \dotplus (0, y_l)$ and, as the sequence $((0, y_1), \ldots, (0, y_n))$ is $\varepsilon^2$-orthogonal, $y_{l+1}(s) = 0$. □

LEMMA 8.14. *Let $\varepsilon \in (0, \frac{1}{6})$ and let $a_1, a_2 \in [-1, 1]$. Let $(\alpha_1, \ldots, \alpha_n)$ be an $\varepsilon^2$-strongly orthogonal sequence of stationary strategy pairs such that $\alpha_k$ is an $\varepsilon$-equilibrium for each $k = 1, \ldots, n$. Assume that for each $k$, $\gamma_k \in [a_1, a_1 + \varepsilon] \times [a_2, a_2 + \varepsilon]$, where $\gamma_k$ is the payoff that corresponds to $\alpha_k$. Let $\alpha = \alpha_1 \dotplus \cdots \dotplus \alpha_n = (x, y)$. Then:*

(a) *We have $a_i - \varepsilon \leq \gamma^i(x, y)$.*
(b) *For each pair $(x', y')$ of stationary strategies, $\gamma^1(x', y) \leq a_1 + 8\varepsilon$ and $\gamma^2(x, y') \leq a_2 + 8\varepsilon$.*

PROOF. Denote $\alpha_k = (x_k, y_k)$. We prove the result only for player 1. We first prove (a). Since $a_1 \leq \gamma_k^1(x_k, y_k)$ for each $1 \leq k \leq n$, it follows from Corollary 8.11 and since $\varepsilon < 1/6$, that $a_1 - \varepsilon \leq a_1 - 6\varepsilon^2 \leq \gamma^1(x, y)$.

We now prove (b). Let $\bar{x}$ be a stationary strategy that maximizes player 1's payoff against $y$: $\gamma^1(\bar{x}, y) = \max_{x'} \gamma^1(x', y)$. Fixing $y$, the game reduces to a Markov decision process and hence such an $\bar{x}$ exists. Moreover, there exists such a strategy $\bar{x}$ that is pure [i.e., $\bar{x}(s) \in \{0, 1\}$ for every $s$] and stops at most once in every branch. Observe that since the sequence $(\alpha_1, \ldots, \alpha_n)$ is $\varepsilon^2$-strongly orthogonal, so is the sequence $((0, y_1), \ldots, (0, y_n))$. Let $\bar{x}_1, \ldots, \bar{x}_k$ be the strategies defined in Lemma 8.13 w.r.t. $\bar{x}$ and $y_1, \ldots, y_n$. Then $\bar{x} = \bar{x}_1 \dotplus \cdots \dotplus \bar{x}_n$ and $(\bar{\alpha}_1, \ldots, \bar{\alpha}_n)$ is $\varepsilon$-orthogonal, where $\bar{\alpha}_k = (\bar{x}_k, y_k)$.

For each $k$, $(x_k, y_k)$ is an $\varepsilon$-equilibrium and, therefore, $\gamma^1(\bar{x}_k, y_k) \leq a_1 + 2\varepsilon$. By Corollary 8.11 and the definition of $\bar{x}$, for every $x'$ we have $\gamma^1(x', y) \leq \gamma^1(\bar{x}, y) \leq a_1 + 2\varepsilon + 6\varepsilon = a_1 + 8\varepsilon$. □

8.5. *Proof of Proposition* 5.5. We now prove Proposition 5.5. Consider the following recursive procedure:

1. *Initialization*: Start with the game $\tilde{T} = T$, the strategy pair $\sigma_0 = (0,0)$ (always continue) and $k = 0$.
2. If there exists a stationary $\varepsilon$-equilibrium in a subgame $T'$ of $\tilde{T}$ with corresponding payoff in $[a_1, a_1 + \varepsilon] \times [a_2, a_2 + \varepsilon]$:
   (a) Set $k = k+1$ and let $\alpha_k = (x_k, y_k)$ be any such $\varepsilon$-equilibrium. Extend $x_k$ and $y_k$ to strategies on $T$ by setting $x_k(s) = y_k(s) = 0$ for every node $s \in S_0 \setminus T'$.
   (b) Set $\sigma_k = \sigma_{k-1} \dotplus \alpha_k$.
   (c) Let $H_k = H_{\varepsilon^2}(\sigma_k)$ be the set of $\varepsilon^2$-heavy nodes of $\sigma_k$ (by Fact 1, $H_{k-1} \subseteq H_k$). Set $\tilde{T} = T_{H_k}$.
   (d) Start stage 2 all over.
3. If, for all subgames $T'$ of $\tilde{T}$, there are no $\varepsilon$-equilibria in $T'$ with corresponding payoff in $[a_1, a_1 + \varepsilon] \times [a_2, a_2 + \varepsilon]$, set $n = k$, $x = x_1 \dotplus \cdots \dotplus x_n$, $y = y_1 \dotplus \cdots \dotplus y_n$ and $D = H_n$.

The idea is to keep adding strongly orthogonal $\varepsilon$-equilibria as long as we can. The procedure continues until there is no $\varepsilon$-equilibrium in any subgame of $\tilde{T}$ with payoffs in $[a_1, a_1 + \varepsilon] \times [a_2, a_2 + \varepsilon]$. The termination of the procedure follows from Lemma 8.4.

The first part of Proposition 5.5 is an immediate consequence of the termination of the procedure. We now prove that $\sigma_n = (x, y)$ satisfies the requirements of the second part. Since $D = H_n$ is the set of $\varepsilon^2$-heavy nodes of $(x, y)$, $\pi(x, y) \geq \varepsilon^2 \times p_D$. For every $1 \leq k \leq n$, $\gamma^i(x_k, y_k) \geq R^i - \varepsilon$, so that $(x_k, y_k)$ is an $\varepsilon$-equilibrium in $T$. Thus $((x_1, y_1), \ldots, (x_n, y_n))$ is an $\varepsilon^2$-strongly orthogonal sequence of stationary $\varepsilon$-equilibria. The remaining claims of Proposition 5.5 follow from Lemma 8.14.

**9. More than two players.** When there are more than two players, it is no longer true that games on a tree admit stationary $\varepsilon$-equilibria. An example of a three-player game where this phenomenon happens was first found by Flesch, Thuijsman and Vrieze (1997). Nevertheless, a consequence of Solan (1999) is that every three-player game on a tree admits a periodic $\varepsilon$-equilibrium, but the period may be long [see Solan (2001)]. We do not know whether this result can be used to generalize Proposition 5.5 for three-player games.

When there are at least four players, the existence of $\varepsilon$-equilibria in stopping games on finite trees is still an open problem, even in the deterministic case; that is, when every node in the tree has at most a single child. For more details, the reader is referred to Solan and Vieille (2001).

**Acknowledgment.** We thank an anonymous referee for helpful comments, which substantially improved the presentation.

School of Mathematical Sciences
Tel Aviv University
Tel Aviv 69978
Israel
e-mail: gawain@post.tau.ac.il

Kellogg School of Management
Northwestern University
and
School of Mathematical Sciences
Tel Aviv University
Tel Aviv 69978
Israel
e-mail: eilons@post.tau.ac.il
e-mail: eilons@kellogg.northwestern.edu
url: www.math.tau.ac.il/eilons
url: www.kellogg.northwestern.edu/faculty/solan/htm/solan.htm